\documentclass{gtart}
\gtart
\usepackage{amssymb,latexsym, amscd}
\usepackage[all]{xy}
\usepackage{graphicx}

 \newtheorem{theorem}{Theorem}[subsection]
 \newtheorem{cor}[theorem]{Corollary}
 \newtheorem{lemma}[theorem]{Lemma}
 \newtheorem{proposition}[theorem]{Proposition}
 \theoremstyle{definition}
 \newtheorem{definition}[theorem]{Definition}
 \theoremstyle{definition}
 \newtheorem{example}[theorem]{Example}
 \theoremstyle{remark}
 \newtheorem{rem}[theorem]{Remark}

\usepackage{latexsym}
\usepackage{amssymb}
\newcommand{\ben}{\begin{equation}}
\newcommand{\een}{\end{equation}}

\newcommand{\integer}{\ensuremath{{\mathbb Z}}}
\newcommand{\naturals}{\ensuremath{{\mathbb N}}}
\newcommand{\real}{\ensuremath{{\mathbb R}}}
\newcommand{\complex}{\ensuremath{{\mathbb C}}}

\newcommand{\U}[1]{\ensuremath{{\mathrm U( #1 )}}}

\newcommand{\Aa}{{\mathcal A}}
\newcommand{\PP}{{\mathcal P}}
\newcommand{\XX}{{\mathcal X}}

\newcommand{\UU}{{\mathcal U}}

\newcommand{\WW}{{\mathcal W}}
\newcommand{\SSS}{{\mathcal S}}
\newcommand{\FF}{{\mathcal F}}

\newcommand{\LL}{\mathcal{L}}

\newcommand{\TT}{\mathcal{T}}

\newcommand{\Xx}{\mathsf{X}}
\newcommand{\Gg}{\mathsf{G}}

\newcommand{\Mm}{\mathsf{M}}
\newcommand{\Ss}{\mathsf{S}}

\newcommand{\Ff}{\mathsf{F}}

\newcommand{\target}{\mathsf{t}}
\newcommand{\source}{\mathsf{s}}
\newcommand{\ident}{\mathsf{e}}
\newcommand{\invers}{\mathsf{i}}
\newcommand{\mult}{\mathsf{m}}

\newcommand{\Loop}{\mathsf{L}}

\newcommand{\twoarrows}{\rightrightarrows}

\newcommand{\To}{\longrightarrow}
\newcommand{\twodownarrows}{\downdownarrows }
\newcommand{\timests}{\: {}_{\target}  \! \times_{\source}}

\newcommand{\Gerbe}{{\texttt{Gb}}}

\newcommand{\coarse}[1] {{#1} / \! \! \sim}
\newcommand{\hyper}{\ensuremath{{\mathbb H}}}

\newcommand{\Diff}{\ensuremath{{\mathrm{Diff}}}}
\newcommand{\VS}{\mathsf{Vector\ Spaces}}

\begin{document}

\title[Holonomy for Gerbes over Orbifolds]
{Holonomy for Gerbes over Orbifolds.}

\author{ Ernesto Lupercio and Bernardo Uribe}
\thanks{The first author was partially supported by the National Science
Foundation}
\address{Department of Mathematics, University of Wisconsin at Madison, Madison, WI 53706}
\address{Max-Planck-Institut f\"{u}r Mathematik, Vivatsgasse
7, D-53111 Bonn, Germany Postal Address: PO.Box: 7280, D-53072
Bonn} \email{ lupercio@math.wisc.edu\footnote{The first author was
partially supported by the National Science Foundation} \\
uribe@mpim-bonn.mpg.de}

\begin{abstract}
In this paper we compute explicit formulas for the holonomy map for
a gerbe with connection over an orbifold. We show that the holonomy
descends to a transgression map in Deligne cohomology. We prove
that this recovers both the inner local systems in Ruan's theory of
twisted orbifold cohomology  \cite{Ruan1} and the local system of
Freed-Hopkins-Teleman in their work in twisted K-theory
\cite{FreedHopkinsTeleman}. In the case in which the orbifold is
simply a manifold we recover previous results of Gaw{\c
e}dzki\cite{Gawedzki} and Brylinski\cite{Brylinski}.
\end{abstract}

 \maketitlepage

\section{Introduction}

A gerbe $\LL$ over a manifold $M$ (or a scheme, if you prefer) has
much in common with a complex line bundle $L$.

A complex line bundle $L$ is classified up to isomorphism class by a
cohomology class $c_1(L) \in H^2(M;\integer)$, its Chern class. By
using the exponential sequence of sheaves
$$ 0 \rightarrow \underline{\integer}
\rightarrow \underline{\complex} \stackrel{\exp}{\rightarrow}
\underline{\complex}^\times \rightarrow 1$$
 we can immediately interpret the Chern class of
$L$ as an element $[g]$ of the cohomology group
$H^1(M;\underline{\complex}^\times)$. In fact a \v{C}ech cocycle
for this class is given by the gluing maps $g_{ij} \colon U_{ij} =
U_i \cap U_j \to \complex^\times$ of the line bundle for a Leray
atlas $(U_i)_{i\in I}$ of the manifold $M$, namely one in which
all open sets and their finite intersections are empty or
contractible.

Moreover, if we put a connection $\nabla$ on the line bundle given
locally by 1-forms $A_i \in \Omega^1(U_i) \otimes \complex$ then
the curvature $F(L,\nabla) = dA \in \Omega^2(M) \otimes \complex$
satisfies the Bianchi identity $$dF=0$$ and therefore it defines a
cohomology class $[F] \in H^2(M; \complex)$. Weil \cite{Weil52}
showed that $[F]$ is the image of $-c_1(L)$ under the map
$H^2(M;\integer) \rightarrow H^2(M;\complex)$.

In any case $c_1(L)$ completely determines the isomorphism class
$[L]$ of $L$ -- we say that $H^2(M; \integer)$ is isomorphic to the
group of isomorphism classes of line bundles over $M$ --. Later in
this paper we will define a cohomology group
$\hyper^2(M;\integer(2)_D^\infty)$ due to Deligne and Brylinski
(and also Cheeger-Simons)  that has the following properties:
   \begin{itemize}
      \item There is a surjective homomorphism $\hyper^2(M;\integer(2)_D^\infty) \rightarrow H^2(M;\integer)$
      \item $\hyper^2(M;\integer(2)_D^\infty)$ classifies isomorphism
      classes $[L,\nabla]$ of line bundles with connection and the
      map above is realized by $[L,\nabla]\mapsto [L]$
   \end{itemize}

 Let us denote by $\LL M$ the space of smooth maps from the
circle $S^1$ to $M$ (with no base point condition -- $\LL M$ is
known as the \emph{free loop space} of $M$). There is a
tautological map $$S^1 \times \LL M \longrightarrow M$$ called the
\emph{evaluation map} sending $(z,\gamma) \mapsto \gamma(z)$. We
can use this map together with the K\"unneth theorem and the fact
that $H^1(S^1;\integer)=\integer$  to get
$$H^2(M; \integer) \rightarrow H^2(S^1 \times \LL M ; \integer) \cong
H^2(\LL M; \integer) \oplus ( H^1(\LL M ; \integer) \otimes
H^1(S^1; \integer) )$$ $$ \stackrel{\cong}{\rightarrow} H^2(\LL M;
\integer) \oplus H^1(\LL M ; \integer) \rightarrow H^1(\LL M ;
\integer)\cong H^0(\LL M ; \underline{\complex}^\times)$$ (where
the next to last map is projection into the second component, and
the last is induced by the exponential sequence). We call the
resulting map $H^2(M; \integer) \rightarrow H^1(\LL M ; \integer)$
the \emph{transgression map.}

The previous discussion can be refined to get a map $$\hyper^2(M;
\integer(2)_D^\infty) \rightarrow H^0(\LL M ;
\underline{\complex}^\times)$$ that has a classical interpretation
in terms of the \emph{holonomy} of $(L,\nabla)$ along a closed
path $\gamma \in \LL M$. To wit, a connection $\nabla$ on $L$
produces a \emph{parallel transport}, that is a linear map
$P_{(L,\nabla)}(\gamma)$ for every path
$\gamma\colon[a,b]\rightarrow M$ of the form
$$P_{(L,\nabla)}(\gamma) \colon L_{\gamma(a)} \rightarrow
L_{\gamma(b)}$$ from the initial fiber to the final fiber.

Let us define the category $\SSS^0(M)$ that we call the \emph{0-th
Segal category of $M$}. Its objects are the points of $M$ and its
morphisms are paths $\gamma\colon[0,2 \pi]\rightarrow M$. We think
of $\gamma$ as an arrow from $\gamma(0)$ to $\gamma( 2 \pi)$. Then
given a line bundle with connection  $(L,\nabla)$ the parallel
transport gives us a functor $$ P_{(L,\nabla)} \colon \SSS^0(M)
\longrightarrow \VS$$ that assigns to the object $x \in M$ the
one-dimensional vector space $L_x$ and to the arrow $\gamma$ the
linear map $P_{(L,\nabla)}(\gamma)$.

In particular, should $\gamma$ be a closed path $\gamma\colon S^1
 \to M$ then the linear isomorphism $P_{(L,\nabla)}(\gamma)$
can be canonically identified with an element of $\complex^\times$,
producing then a map
$$\LL M \longrightarrow \complex^\times$$ and hence and element in
$H^0(\LL M ; \underline{\complex}^\times)$.

Let us consider now a gerbe $\LL$ over $M$. We will define gerbes
later in the paper, but for now we list some of their properties.
\begin{itemize}
   \item The group of isomorphism classes of gerbes $\Gerbe(M)$ on
   $M$ is isomorphic to $H^3(M;\integer)$.
   \item The isomorphism $\Gerbe(M)\rightarrow H^3(M;\integer)$
   is realized by the \emph{Dixmier-Douady} characteristic
   class $dd(\LL)\in H^3(M;\integer)$ of the gerbe $\LL$.
   \item We can place connections $\Xi$ (also known as connective
   structures) on gerbes.
   \item The curvature of a connection over a gerbe $\LL$ on $M$ is
   a closed 3-form $G\in\Omega^3(M) \otimes \complex$.
   \item The de Rham cohomology class $[G] \in H^3(M;\real)$ is the
   real image of $dd(\LL)$.
   \item The group of isomorphism classes of gerbes with
   connections over $M$ is isomorphic to a \emph{Deligne
   cohomology} group $\hyper^3(M;\integer(3)_D^\infty)$.
   \item The holonomy of a gerbe $\LL$ with connection $\Xi$ is a
   complex line bundle $L$ with connection $\nabla$ on the free loop space $\LL M$.
   \item The holonomy $(\LL,\Xi) \mapsto (L,\nabla)$ realizes
   a transgression map
   $$\hyper^3(M;\integer(3)_D^\infty) \rightarrow \hyper^2(\LL M;\integer(2)_D^\infty)$$
   \item A pair $(\LL, \Xi)$ induces a parallel transport functor
   $$P_{(\LL,\Xi)}\colon \SSS^1(M)\longrightarrow \VS$$ from the
   first Segal category of $M$ whose objects are maps from compact
   closed one-dimensional oriented manifolds to $M$, and whose morphisms
   are maps from compact 2-dimensional manifolds to $M$ forming cobordisms
   between two objects \cite{Murray02,Segal02}. For instance, in
   the picture below we have two maps $\gamma_i \colon S^1
   \rightarrow M$ ($i=1,2$) and a map $\Sigma \colon F \rightarrow
   M$ from a 2-dimensional manifold $F$ into $M$. Such a
   configuration would produce a linear isomorphism
   $$P_{(\LL,\Xi)}(\Sigma) \colon L_{\gamma_1} \rightarrow
   L_{\gamma_2}$$

\begin{eqnarray}
\includegraphics[height=1.0in]{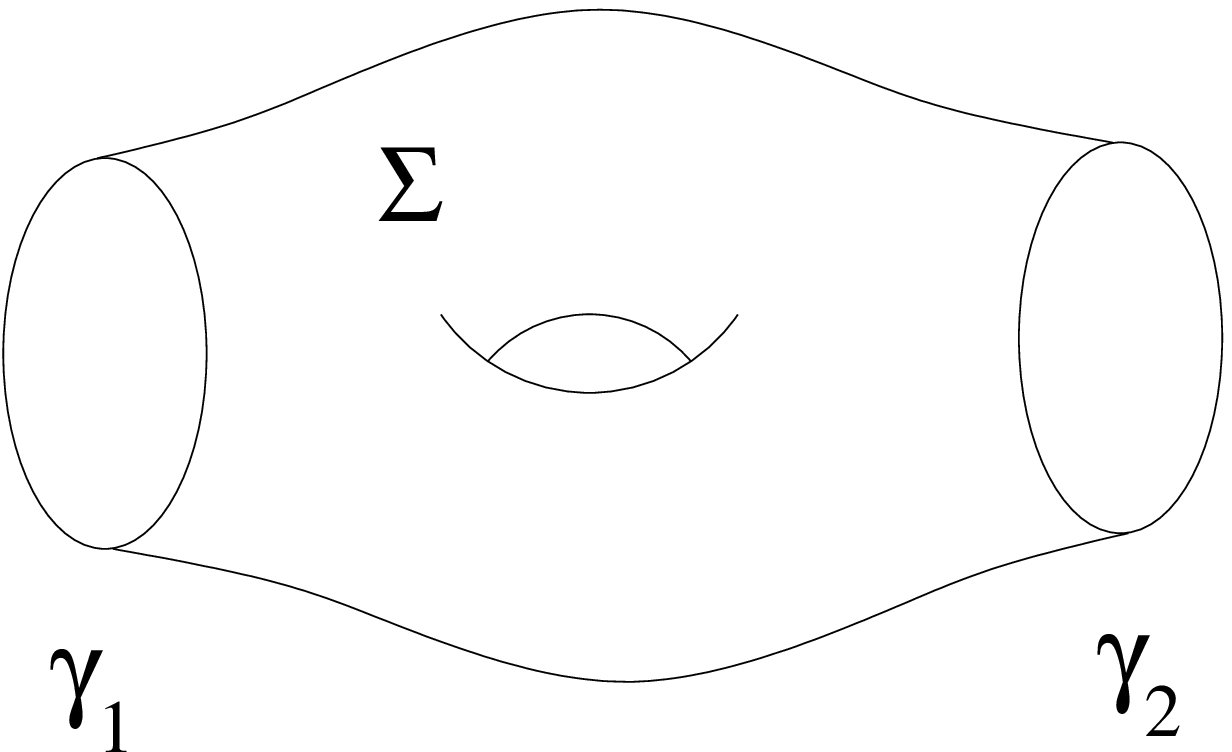} \label{graph string}
\end{eqnarray}

 Such a functor is closely related to  a \emph{String Connection} in the
   terminology of Segal \cite{Segal02}, Stolz, and Teichner \cite{StolzTeichner1}.
   More specifically if $\gamma
   \colon S^1 \to M$ is an object of $\SSS^1(M)$,
   then we have $P_{(\LL,\Xi)}(\gamma) = L_\gamma$, where $L$ is the line
   bundle over $\LL M$ mentioned in the last paragraph.
\end{itemize}

The purpose of this paper is to generalize the previous picture to
the case in which instead of considering a manifold $M$ we
consider a smooth \emph{orbifold} or \emph{Deligne-Mumford stack}
$\XX$. This new case involves many new features and links together
several interesting structures that have appeared in geometry and
topology recently. We will briefly describe now the contents of
this paper.

In Section \ref{section Deligne Cohomology for Groupoids}
 we set our notations and terminology for the theory of
groupoids. We will use groupoids as models for our orbifolds --
they will be our basic tool. We recall that a groupoid $\Gg$ is a
category in which all morphisms have inverses and by an orbifold groupoid
we mean a proper, smooth, \'{e}tale groupoid. In this paper
whenever we write groupoid a {\emph{smooth, \'{e}tale groupoid}} 
is to be understood.

In Section \ref{sectionloopgroupoid} we deal with the issue of
defining the "free loop space" of an orbifold. Since an orbifold
is already more than a space, the answer is itself an infinite
dimensional orbifold that we call \emph{the loop orbifold}. Our
model for the loop orbifold will be a groupoid -- the \emph{loop
groupoid} $\Loop \Gg$. In particular when the orbifold happens to
be a manifold then the loop orbifold is simply the free loop space
of the manifold.

In Section \ref{subsection Sheaves and cohomology} we explain
 sheaf cohomology theory for groupoids. Then in section
\ref{subsection Deligne Cohomology} we use this theory to define
Deligne cohomology for groupoids and explain its relation to the
theory of $n$-gerbes with connective structure.

In Section \ref{section Holonomy} we prove the following theorem.

\begin{theorem}
There is a natural transgression map (holonomy)
$$\tau_1 :\breve{C}^1(\Gg,\complex^\times(2)_\Gg) \To \breve{C}^0(\Loop \Gg,\complex^\times_{\Loop \Gg})$$
that associates to every line bundle with connection over $\Gg$
its holonomy. Here $\complex^\times_{\Loop \Gg}$ is the sheaf of
$\complex^\times$ valued functions on the loop groupoid. This map descends to cohomology
$$ \hyper^1(\Gg;\integer(1)_D^\infty)
 \To H^0\left(\Loop \Gg, \complex^\times_{\Loop \Gg}\right).$$
\end{theorem}

In Section \ref{section Line bundle over the loop groupoid}
 we go ahead and \emph{define} the Deligne cohomology
for the loop groupoid to then prove the following theorem.

\begin{theorem}
There is a natural holonomy homomorphism
$$\tau_2 :\breve{C}^2(\Gg,\complex^\times(3)_\Gg) \To \breve{C}^1(\Loop \Gg,\complex^\times(2)_{\Loop \Gg})$$
from the group of gerbes with connection over $\Gg$ to the group of
line bundles with connection over the loop groupoid. Moreover this
  holonomy map commutes with the coboundary operator
and therefore induces a map in orbifold Deligne cohomology
$$\hyper^2(\Gg;\integer(2)_D^\infty)
\To \hyper^1(\Gg;\integer(1)_D^\infty).$$
\end{theorem}

We prove in fact a little bit more. For what we really construct is
a functor  $$\Loop \Gg \rightarrow \VS$$ given by parallel
transport along arrows of the loop groupoid for the gerbe
connection. While in the case of a manifold $M$ this is only the
portion of the string connection that associates the vector space
to the objects of $\SSS^1(M)$, in the case of an orbifold groupoid
$\Gg$ we already have arrow assignments. Since $\Loop \Gg
\hookrightarrow \SSS^1(\Gg)$ is an inclusion of categories we
think of the functor constructed here as a genus-zero ghost part
of the string connection. We will return to the construction of
the full string connection for $\SSS^1(\Gg)$ in a future paper.

Furthermore in Section \ref{subsection Manifolds} we show that
when the orbifold in question is simply a manifold then we recover
the results of Gaw{\c e}dzki \cite{Gawedzki} and Brylinski
\cite{Brylinski}.

In Section \ref{subsection Global quotients} we study the
particular case in which the orbifold $\Xx$ is actually a global
quotient $[M/G]$, and the gerbe in question is globally defined.

Then in Section \ref{subsection Discrete torsion} we compute
explicitly the case when the gerbe in question comes from discrete
torsion in a global quotient orbifold $[M/G]$.

In Section \ref{section Localization at the fixed points}
 we pursue the subject of localization. One of the main
results of \cite{LupercioUribe2} is the following theorem

\begin{theorem} The fixed suborbifold of $\Loop \Gg$ under the
natural $S^1$-action (rotating the loops) is $$ \wedge \Gg = (\Loop
\Gg)^{S^1}$$
\end{theorem}

Here $\wedge \Gg$ is the so-called \emph{twisted sector orbifold}
or \emph{inertia groupoid} of $\Gg$. In the case of a manifold
$\wedge M = M$.

We prove then the following theorem for $\Gg$ an orbifold groupoid

\begin{theorem} The restriction of the holonomy of a gerbe with
connection over $\Gg$ (that is a line bundle with connection over
$\Loop \Gg$) is an inner local system on $\wedge \Gg$.
\end{theorem}

Inner local systems were discovered by Ruan \cite{Ruan1} for
completely different reasons. As it happens this is too the local
system used by Freed, Hopkins and Teleman
\cite{FreedHopkinsTeleman} in their work on twisted $K$-theory
whenever the action of the Lie group is almost free, namely that it has
only finite stabilizers.

Finally in Section 6 we discuss how to generalize the previous
theory to $n$-gerbes with connection. The corresponding holonomy
formula send $(n+1)$-gerbes to $n$-gerbes.

What we prove in this paper is actually a bit stronger than the
statements of the previous theorems. We give \emph{explicit
formulas} for the holonomy maps, and then show that it descends to
Deligne cohomology. Our motivation to do this is that in physics
all the objects we have discussed have interesting interpretations
and explicit formulas are necessary for computations
\cite{FreedWitten}. For example we have explained elsewhere  that
in orbifold string theory (and conformal field theory) both the
$B$-field and discrete torsion can be suitably interpreted in
terms of gerbes and Deligne cohomology over orbifolds
\cite{LupercioUribe3}. For related statements and work in the
physics literature we refer the reader to \cite{Sharpe} and the
references therein.

{\textbf{ Acknowledgments. } } We would like to thank
conversations with A. Adem, M. Ando, A. Carey, D. Freed, M.
Hopkins, P. Lima-Filho, J. Mickelsson, Y. Ruan, A. Schilling, G.
Segal, S. Stolz, C. Teleman and A. Waldron. Both authors would
like to thank the Erwin Schr\"{o}dinger International Institute for
Mathematical Physics where part of the work for this paper took
place. The second author would like to thank the hospitality of
the Max Planck Institut f\"{u}r Mathematik in Bonn where this
paper took its final form.
Finally we would like to specially thank D. Freed for valuable
correspondence and comments regarding the first draft of this paper.

\section{Deligne Cohomology for Groupoids} \label{section Deligne Cohomology for Groupoids}

When we say a groupoid we mean a (small) category $\Gg$ so that
the set of its objects $\Gg_0$ and the set of its arrows $\Gg_1$
are both manifolds, and every arrow has an inverse.

We will represent orbifolds by groupoids. It is useful to consider
the following two examples as motivation.

\begin{example} Let $G$ be a finite group. Consider the orbifold $[M/G]$
obtained from a $G$-manifold $M$ (we use the brackets to
differentiate the orbifold $[M/G]$ from the quotient space -- or
coarse moduli of orbits -- $M/G$ ). Then we will associate to it
the groupoid $\Xx$ with morphisms $\Xx_1 = M \times G$ and objects
$\Xx_0 =M$. The arrow $(m,g)$ takes the object $m$ to the object
$mg$. We will often write the groupoid $\Xx$ by the symbol $M\times
G \rightrightarrows M$.
\end{example}

\begin{example}\label{manifoldgroupoid} Consider a manifold $M$ with an atlas $\UU =
(U_i)_i$. We will associate to $(M,\UU)$ the groupoid
$\Mm_\UU$ whose objects $\Mm_0 = \{ (x,i) : x \in U_i \} =
\coprod_i U_i$ and whose arrows $\Mm_1 = \{ (x,i,j) : x \in U_{ij}
= U_i \cap U_j \} = \coprod_{(i,j)} U_{ij}$. The arrow $(x,i,j)$
takes the object $(x,i)$ to the object $(x,j)$.
\end{example}

In the case of a general orbifold there is a groupoid that
represents it that is a sort of hybrid of the previous two
examples.

The groupoids $\Gg$ we will be concerned with will be \'{etale}
and smooth, this means that all the structure maps
 $$\xymatrix{
         \Gg_1 \timests \Gg_1 \ar[r]^{\mult} & \Gg_1 \ar[r]^{\invers} &
         \Gg_1 \ar@<.5ex>[r]^{\source} \ar@<-.5ex>[r]_{\target} & \Gg_0 \ar[r]^{\ident} & \Gg_1
         }$$
( $\mult=$ composition, $\invers=$ inverse, $\ident=$ identity,
$\source=$ source, $\target=$ target) are local diffeomorphisms
with $\Gg_i$ manifolds. When the anchor map $(\source, \target) :
\Gg_1 \to \Gg_0 \times \Gg_0$ is proper the groupoid will
represent an orbifold groupoid.
\begin{definition}
By an orbifold groupoid we mean smooth,
\'{e}tale, proper groupoid.
\end{definition}
 On this paper we are only concerned
with smooth \'{e}tale groupoids. When we think of a
groupoid, implicitly what we are considering, is the Morita
equivalence class where the groupoid belongs (see \cite{Moerdijk2002} for
the definition of Morita equivalence). But in order to make
calculations explicit, or to use {\v C}ech cohomology, we will
make use of a special representative of the Morita class. We
require this groupoid to be built out of a disjoint union of
contractible sets as follows.
\begin{definition} \label{Leray}
A groupoid $\Gg$ is called {\emph {Leray}} if $\Gg_i$ is
diffeomorphic to a disjoint union of contractible open sets for
all $i \in \naturals$.
\end{definition}
The existence of such Leray groupoid representative for every
orbifold is proven by Moerdijk and Pronk  \cite[Cor.
1.2.5]{MoerdijkPronk1}.

Here we are concerned with the geometry of the groupoid and we
will give very explicit geometric descriptions of the objects in
study. The algebraic topology  of  groupoids has been studied by
several authors \cite{CrainicMoerdijk, LupercioUribe1,
LupercioUribe5, Moerdijk2002, MoerdijkPronk1}  and having both the
geometric and the topological approaches is very useful.

Here we should introduce another very important structure
associated to a groupoid called inertia groupoid
\begin{definition} \label{def.inertiagroupoid}
The inertia groupoid $\wedge \Gg$ is defined by:
\begin{itemize}
\item Objects $(\wedge \Gg)_0$: Elements $v \in \Gg_1$ such that $s(v) = t(v)$.
\item Morphisms $(\wedge \Gg)_1$: For  $v,w \in (\wedge \Gg)_0$ an
arrow $v \stackrel{\alpha}{\to} w$ is an element $\alpha \in
\Gg_1$ such that $v \cdot \alpha = \alpha \cdot w$
        $$
         \xymatrix{
         \circ \ar@(ul,dl)[]|{v} \ar@/^/[rr]|{\alpha}
         &&\circ \ar@(dr,ur)[]|{w^{-1}} \ar@/^/[ll]|{\alpha^{-1}}
         }$$
\end{itemize}
\end{definition}
It is known that the inertia groupoid in the case of an orbifold
matches with what is commonly known in the literature by twisted
sectors (see \cite{Ruan1, LupercioUribe2}),
 thus this is
a natural way to define them. And we will see in the next section
that the inertia groupoid arises naturally as the constant loops
of the loop groupoid.

\subsection{Loop Groupoid} \label{sectionloopgroupoid}

In the last years it has become increasingly evident that the free
loop space of a manifold $\LL M =Map(S^1,M)$ is a very important
concept, providing a sort of natural thickening of $M$, it is at
first not clear how one must define the free loop space of an
orbifold. It is apparent that simply considering free loops in the
quotient space $\coarse{\Gg}$ forgets all the orbifold structure.
In fact the correct notion that we call the \emph{loop orbifold}
or \emph{loop groupoid} is itself an orbifold, albeit an infinite
dimensional one.

We have defined in \cite{LupercioUribe2} the loop groupoid to be a
category whose objects are $Hom(S^1,\Gg)$ in the category of
groupoids (its objects were known to Mr{\v{c}}un \cite{Mrcun} and to
Bridson-Haefliger \cite{BridsonHaefliger}).

\begin{example} \label{example loop groupoid}
Consider again the orbifold $X = [M/G]$ represented by the
groupoid $M\times G \rightrightarrows M$ that we denote by $\Xx$,
as at the beginning of the last section. Consider the groupoid
$\Loop \Xx$ whose objects are all pairs $(\phi,g)$, $\phi: [0,1]
\to M$, $g \in G$ where we have $\phi(0) g = \phi(1)$. Let $G$ act
on the paths in the natural way, i.e. $\{\phi \cdot k\} (t) =
\phi(t) k$ with $\{\phi\cdot k\}(0) k^{-1}gk =\{\phi \cdot
k\}(1)$. We declare the arrows of $\Loop \Xx$ to be the triples
$(\phi,g,k)$ so that $\phi: [0,1] \to M$, $\phi(0) g = \phi(1)$
and $k \in G$. The arrow $\Lambda = (\phi, g, k)$ in $\Loop \Xx$
sends the path $(\phi, g)$ to $(\phi \cdot k, k^{-1}gk)$. We call
$\Loop \Xx$ \emph{the loop groupoid}.
\end{example}

To do this in full generality we must face the following
difficulty. Suppose that we first assign a groupoid to $S^1$ and
consider the groupoid morphisms to $\Gg$. They will  certainly be a portion
of the desired loop groupoid, but unfortunately we may be missing
elements on it that will only become apparent by choosing a finer
groupoid representation of the circle. Therefore we need to
consider all the Morita equivalent groupoids representing the
circle (this amounts to take finer and finer covers of the
circle). This is explained in detail in \cite{LupercioUribe2}. The
following formalism (that is unfortunately a bit technical) solves
this difficulty.

For a finite set $\{q_1, \dots, q_n, q_0\} \subset (0,1]$ with
$q_1 < \cdots < q_n <q_0$ and $\epsilon > 0$ sufficiently small we
associate a unique cover of the circle given by the sets
$V_i^0:=(q_i - \epsilon, q_{i+1} + \epsilon)$ and the exponential
map $e^{2 \pi i t}$. This cover induces an {\it admissible cover}
$W$ on the real numbers $\real$ that consist of the sets
$V_i^k:=(q_i+k - \epsilon, q_{i+1} + k + \epsilon)$ and
$V_{0}^k:=(q_0+k-1 - \epsilon, q_{1} + k+ \epsilon)$ for $k \in
\integer$ and $1 \leq i \leq n$. We call $\real^W$ the groupoid
associated to this cover, i.e.
$$\real^W_1 : = \bigsqcup_{i,j,k,l} V_{i,j}^{k,l} \ \ \ \ \  \real^W_0 := \bigsqcup_{i,k}V_i^k$$
where $V_{i,j}^{k,l}:= V_i^k \cap V_j^l$; and the $\epsilon$ is
small enough so that all the double intersections are of the form
$(q_i +k - \epsilon, q_i +k + \epsilon)$ or empty. We can now define the
natural action of $\integer$ in $\real^W_1$
\begin{eqnarray*}
\real^W_1 \times \integer & \to & \real^W_1\\
((x, V_{i,j}^{k,l}),m) & \to & (x+m, V_{i,j}^{k+m,l+m})
\end{eqnarray*}
with $x \in V_{i,j}^{k,l}$ and $x + m \in V_{i,j}^{k+m,l+m}$.
\begin{definition}
Let $\Ss^1_W$ be the groupoid
$$\begin{array}{c}
\real^W_1 \times \integer\\
\twodownarrows \\
\real^W_0
\end{array}$$
with maps
$$\source\left((x,V_{i,j}^{k,l}),m\right)=(x,V_i^k) \ \ \ \ \ \ \
\target \left((x,V_{i,j}^{k,l}),m\right)=(x+m,V_{j}^{l+m})$$
$$\ident(x,V_i^k) = \left((x,V_{i,i}^{k,k}),0\right) \ \ \ \ \ \ \
\invers\left((x,V_{i,j}^{k,l}),m\right)=\left((x+m,V_{j,i}^{l+m,k+m}),-m\right)$$
$$\mult \left[ \left((x,V_{i,j}^{k,l},m\right),\left((x+m,V_{i,j}^{k+m,l+m}),n\right) \right]=
\left((x,V_{i,j}^{k,l}),m+n\right).$$
\end{definition}
The groupoid $\Ss^1_W$ is Morita equivalent to unit groupoid $S^1
\twoarrows S^1$ (all arrows are identities). If $W'$ is a
refinement of $W$, then there is a unique Morita morphism
$\rho^{W'}_W : \Ss^1_{W'} \to \Ss^1_W$.
\begin{definition}
For $\Gg$ a topological  groupoid and an open cover $W$ of the
circle, the loop groupoid $\Loop\Gg(W)$ associated to $\Gg$ and
the open cover $W$ will be defined by the following data:
\begin{itemize}
\item $\Loop\Gg(W)_0$ the objects:
Morphisms of groupoids $\Ss_W^1 \to \Gg$
\item $\Loop\Gg(W)_1$ the morphisms: For two elements in $\Loop\Gg(W)_0$, say
$\Psi, \Phi : \Ss^1_W \to \Gg$ , a morphism (arrow) from $\Psi$ to
$\Phi$ is a map $\Lambda : \real^W_1 \times \integer\to \Gg_1$
that makes the following diagram commute
      $$
       \xymatrix{
        \real^W_1 \times \integer \ar[r]^\Lambda \ar[d]_{\source \times \target}
       & \Gg_1 \ar[d]^{\source \times \target}\\
       \real^W_0 \times \real^W_0 \ar[r]_{(\Psi_0 , \Phi_0)}
  & \Gg_0 \times \Gg_0 }
      $$
and such that for $r \in \real^W_1 \times \integer$
\begin{eqnarray}  \label{eqnbetweenmorphisms}
\Lambda(r)=\Psi_1(r) \cdot \Lambda(\ident \target(r))=
\Lambda(\ident \source(r)) \cdot \Phi_1(r).
\end{eqnarray}
\end{itemize}
\end{definition}
The composition of morphisms is defined pointwise, in other words,
for $\Lambda$ and $\Omega$ with
       $$\xymatrix{
       \Psi \ar@/^/[r]^\Lambda  & \Phi \ar@/^/[r]^\Omega & \Gamma}
       $$
we set $$\Omega \circ \Lambda(\ident \source(r)) := \Lambda(\ident
\source(r)) \cdot \Omega(\ident \source(r))$$ and
 $$\Omega \circ \Lambda(r):= \Omega \circ \Lambda (\ident \source(r)) \cdot \Gamma(r)=\Psi(r) \cdot
\Omega \circ \Lambda (\ident \target(r))$$ The previous properties
imply that an arrow $\Lambda$ determines its source $\Psi$ and its
target $\Phi$. We consider $\Loop\Gg(W)_1$ as a subspace of the
space of smooth maps $Map(\real^W_1, \Gg_1) \times Map(\real^W_0,
\Gg_1)$; in this way $\Loop\Gg(W)_1$ and $\Loop\Gg(W)_0$ inherit
the compact-open topology, making the groupoid $\Loop\Gg(W)$ into a
topological one. For two admissible covers $W_1,W_2$ associated to
$\{q_1, \dots, q_n, q_0\}$ and $\{q'_1, \dots, q'_{n'}, q'_0\}$
respectively, there is always a common refinement. We could take
$W$ associated to $\{q_1, \dots, q_n, q_0\} \cup \{q'_1, \dots,
q'_{n'}, q'_0\}$ and the natural morphisms $\rho^W_{W_i} : \Ss^1_W
\to \Ss^1_{W_i}$. These induce natural monomorphisms of
topological groupoids $\Loop \Gg(W_i) \hookrightarrow \Loop \Gg(W)$
   We want
that two objects $\Psi_i : \Ss^1_{W_i} \to \Gg_1$ ($i = 1,2$) to
be equivalent if the following square is commutative
 $$\xymatrix{
      \Ss^1_{W} \ar[r]^{\rho^{W}_{W_1}} \ar[d]_{\rho^{W}_{W_2}} & \Ss^1_{W_1} \ar[d]^{\Psi_1}\\
      \Ss^1_{W_2} \ar[r]_{\Psi_2} & \Gg. }
      $$
so we define
\begin{definition}
The loop groupoid $\Loop\Gg$ of the groupoid $\Gg$ is defined as
the monotone union (colimit) of the groupoids $\Loop\Gg(W)$ where
$W$ runs over the set $\WW$ of admissible covers
$$\Loop\Gg := \lim_{\overrightarrow{W\in\WW}} \Loop\Gg(W).$$
\end{definition}

In this way the loop groupoid is naturally endowed with a
topology, becoming a topological groupoid. Now we list some facts
about the loop groupoid that can be found in \cite{LupercioUribe2}.
\begin{itemize}
\item For $\Gg$ \'{e}tale and proper, then $\Loop \Gg$ is also \'{e}tale and has finite isotropy.
\item A morphism of groupoids $\Ff \to \Gg$ induces naturally another one at the level of loops $\Loop \Ff \to \Loop \Gg$.
\item If the morphism  $\Ff \to \Gg$ is Morita, then $\Loop \Ff \to \Loop \Gg$ is also Morita.
\item The loop groupoid $\Loop \Gg$ can be endowed with an action of $\real$ in natural way, i.e. shifting the morphisms by $t \in \real$.
The fixed point set groupoid $\Loop \Gg ^\real$ under this action
is Morita equivalent to the inertia groupoid $\wedge \Gg$ (Def.
\ref{def.inertiagroupoid}).
\end{itemize}

\subsubsection{The tangent loop groupoid}

The loop groupoid is endowed with a natural tangent groupoid in the
same way the groupoid $\Gg$ is endowed with its tangent groupoid
$T \Gg$
 defined as $T \Gg_1 \twoarrows T \Gg_0$ with the induced
structure maps (clearly $T\Gg$  is also  smooth and \'{e}tale).
\begin{definition} \label{tangent loop groupoid}
 For $\Psi$
and object of $\Loop \Gg(W)$, the tangent space $T_{\Psi} \Loop
\Gg(W)$ will consist of all morphisms $\xi : \Ss^1_W \to T \Gg$
such that $p \circ \xi = \Psi$ where $p: T \Gg \to \Gg$ is the
natural projection morphism; these will be the objects of $T \Loop
\Gg(W)$. The morphisms of $T \Loop \Gg(W)$ are the natural ones,
namely for an arrow $\Lambda: \Psi \to \Phi$ with $\xi \in T_\Psi
\Loop \Gg(W)$ and $\zeta \in T_\Phi \Loop \Gg(W)$,
 a tangent morphism in $T_\Lambda \Loop \Gg(W)$ between
 $\xi$ and $\zeta$ is a map $\nu:\real^1_W \times \integer \to T \Gg_1$
that makes the following diagram commute
 $$
       \xymatrix{
        \real^W_1 \times \integer \ar[r]_\nu \ar[d]_{\source \times \target} \ar@/^1.5pc/[rr]^{\Lambda}
       & T\Gg_1 \ar[r]_p \ar[d]^{\source \times \target} & \Gg_1 \ar[d]^{\source \times \target}\\
       \real^W_0 \times \real^W_0 \ar[r]^{(\xi,\zeta)} \ar@/_1.5pc/[rr]_{(\Psi_0 , \Phi_0)} & T \Gg_0 \times T \Gg_0 \ar[r]^{(p,p)}
  & \Gg_0 \times \Gg_0 }
      $$
\end{definition}

Taking the inverse limit over the admissible  covers $W$ of $T \Loop \Gg (W)$ we obtain $T\Loop \Gg$.

\subsection{Sheaves and cohomology} \label{subsection Sheaves and cohomology}

All the properties of sheaves and cohomologies of topological
spaces can be extended for the case of smooth \'{e}tale groupoids. This is
done in \cite{Haefliger,CrainicMoerdijk}. Let us briefly summarize
the theory. A $\Gg$-sheaf $\FF$ is a sheaf over $\Gg_0$, namely a
topological space with a projection $p: \FF \to \Gg_0$ which is a
local homeomorphism on which $\Gg_1$ acts continuously. This means
that for $a \in \FF_x=p^{-1}(x)$ and $g \in \Gg_1$ with
$\source(g)=x$, there is an element $ag$ in $\FF_{t(g)}$ depending
continuously on $g$ and $a$. The action is a map $\FF \: {}_{p} \!
\times_s \Gg_1 \to \FF$. For $\FF$ a $\Gg$-sheaf, a section
$\sigma : \Gg_0 \to \FF$ is called invariant if $\sigma(x) g
=\sigma(y)$ for any arrow $x \stackrel{g}{\to} y$.
$\Gamma_{inv}(\Gg,\FF)$ is the set of invariant sections and it
will be an abelian group if $\FF$ is an abelian sheaf. For an
abelian $\Gg$ sheaf $\FF$, the {\it cohomology} groups
$H^n(\Gg,\FF)$ are defined as the cohomology groups of the complex:
$$\Gamma_{inv}(\Gg, \TT^0) \to \Gamma_{inv}(\Gg,\TT^1) \to \cdots$$
where $\FF \to \TT^0 \to \TT^1 \to \cdots$ is a resolution of
$\FF$ by injective $\Gg$-sheaves. When the abelian sheaf $\FF$ is
locally constant (for example $\FF =\integer$) is a result of
Moerdijk \cite{Moerdijk98} that
$$H^*(\Gg,\FF) \cong H^*(B\Gg, \FF)$$
where the left hand side is sheaf cohomology and the right hand
side is simplicial cohomology. There is a {\it basic spectral
sequence} associated to this cohomology. Pulling back $\FF$ along
\begin{eqnarray}
\epsilon_n : \Gg_n \to \Gg_0 \label{pullbacksheaves}
\end{eqnarray}
 $$\epsilon_n(g_1, \dots, g_n) = \target(g_n)$$ it
induces a sheaf $\epsilon_n^*\FF$ on $\Gg_n$ (where the $\Gg$
action on $\Gg_n$ is the natural one, i.e. $(g_1, \dots, g_n)h=
(g_1, \dots, g_nh)$; $\Gg_n$ becomes in this  way a $\Gg$-sheaf)
such that for fixed $q$ the groups $H^q(\Gg_p , \epsilon_p^* \FF)$
form a cosimplicial abelian group, inducing a spectral sequence:
$$H^pH^q(\Gg_\bullet, \FF) \Rightarrow H^{p+q}(\Gg,\FF)$$
So if $0 \to \FF \to \FF^0 \to \FF^1 \to \cdots$ is a resolution
of $\Gg$-sheaves with the property that $\epsilon_p^* \FF^q$ is an
acyclic sheaf on $\Gg_p$, then $H^*(\Gg,\FF)$ can be computed from
the double complex $\Gamma(\Gg_p, \epsilon_p^* \FF^q)$. We
conclude by introducing the algebraic gadget that will allow us to
define Deligne cohomology. Let $\FF^\bullet$ be a cochain complex
of abelian sheaves, then the {\it hypercohomology groups}
$\hyper^n(\Gg,\FF)$ are defined  as the cohomology groups of the
double complex $\Gamma_{inv}(\Gg, \TT^\bullet)$ where $\FF^\bullet
\to \TT^\bullet$ is a quasi-isomorphism into a cochain complex of
injectives.

\subsection{Deligne Cohomology} \label{subsection Deligne Cohomology}
In what follows we will define the smooth Deligne cohomology of a
smooth \'{e}tale groupoid; we will extend the results of
Brylinski\cite{Brylinski} to groupoids and will follow very
closely the description given in there. We will assume that $\Gg$
is Leray (Def. \ref{Leray}) and that the set of objects $\Gg_0$ is
of bounded cohomological dimension. Deligne cohomology is related
to the De Rham cohomology. We will consider the De Rham complex of
sheaves and we will truncate it at level $p$; what interests us is
the degree $p$ hypercohomology classes of this complex. To be more
specific,  let $\integer(p):= (2 \pi \sqrt{-1})^p \cdot \integer$
be the cyclic subgroup of $\complex$ and $\Aa^p_{\Gg, \complex}$ the $\Gg$-sheaf of
complex-valued differential $p$-forms; as $\Gg$ is a smooth
\'{e}tale groupoid the maps $\source$ and $\target$ are local
diffeomorphisms,
 then the action  of $\Gg$ into the sheaf over $\Gg_0$
of complex-valued differential $p$-forms is the natural one given
by the pull back of the corresponding diffeomorphism. Let
$\integer(p)_\Gg$ be the constant $\integer(p)$-valued
$\Gg$-sheaf, and $i : \integer(p)_\Gg \to \Aa^0_{\Gg,\complex}$
the inclusion of constant into smooth functions.
\begin{definition}
Let $\Gg$ be a smooth \'{e}tale groupoid. The {\bf smooth Deligne
complex} $\integer(p)_D^\infty$ is the complex of $\Gg$-sheaves:
$$\integer(p)_\Gg \stackrel{i}{\To} \Aa^0_{\Gg,\complex} \stackrel{d}{\To}
\Aa^1_{\Gg,\complex} \stackrel{d}{\To} \cdots \stackrel{d}{\To}
\Aa^{p-1}_{\Gg,\complex}$$ The hypercohomology groups
$\hyper^q(\Gg,\integer(p)_D^\infty)$ are called the {\bf smooth
Deligne cohomology} of $\Gg$.
\end{definition}

In order to make the explanations clearer, where are going to work
with a quasi-isomorphic complex of sheaves to the Deligne one,
which is a bit simpler.
\begin{definition} Let $\complex^\times(p)_\Gg$ be the following complex of sheaves:
$$\complex^\times_\Gg \stackrel{d \log}{\To} \Aa^1_{\Gg,\complex} \stackrel{d}{\To}
\cdots \stackrel{d}{\To} \Aa^{p-1}_{\Gg,\complex}$$
\end{definition}
There is a quasi-isomorphism between the  complexes $(2\pi
\sqrt{-1})^{-p+1}\cdot \integer(p)_D^\infty$ and
$\complex^\times(p)_\Gg[-1]$ (this fact is explained in Brylinski
\cite{Brylinski} page 216)
    $$\xymatrix{
    (2\pi \sqrt{-1})^{-p+1}\cdot\integer(p)_\Gg \ar[r] &
     \complex_\Gg \ar[r]^{d} \ar[d]_\exp &\Aa^1_{\Gg,\complex} \ar[r]^{d} \ar[d] &
\cdots \ar[r]^{d} & \Aa^{p-1}_{\Gg,\complex} \ar[d]\\
    &    \complex^\times_\Gg \ar[r]^{d \log} &\Aa^1_{\Gg,\complex} \ar[r]^{d} &
\cdots \ar[r]^{d} &\Aa^{p-1}_{\Gg,\complex}
    }$$
hence there is an isomorphism of hypercohomologies:
\begin{eqnarray} \label{isodelignecohomology}
\hyper^{q-1}(\Gg, \complex^\times(p)_\Gg) \cong (2\pi
\sqrt{-1})^{-p+1} \cdot \hyper^q(\Gg, \integer(p)_D^\infty)
\end{eqnarray}

Now let $\Gg$ be a Leray groupoid. We are going to define the
{\v{C}}ech double complex associated to the $\Gg$-sheaf complex
$\complex^\times(p)_\Gg$. Consider the space
$$C^{k,l}=\breve{C}(\Gg_k, \Aa^l_{\Gg,\complex}):=\Gamma(\Gg_k, \epsilon_k^*\Aa^l_{\Gg,\complex})$$
 of global sections of the sheaf $\epsilon_k^*\Aa^l_{\Gg,\complex}$ over
$\Gg_k$ as in (\ref{pullbacksheaves}). The vertical differential
$C^{k,l} \to C^{k,l+1}$ is given by the maps of the complex
$\complex^\times(p)_\Gg$
 and the horizontal differential $C^{k,l} \to C^{k+1,l}$ is
obtained by $\delta = \sum(-1)^i\delta_i$ where for $\sigma \in
\Gamma(\Gg_k, \epsilon_k^*\Aa^l_{\Gg,\complex})$
$$ (\delta_i \sigma)(g_1,\dots,g_{k+1}) = \left\{
\begin{array}{cc}
\sigma(g_1, \dots, g_k)\cdot g_{k+1} & \hbox{for $i=k+1$}\\
\sigma(g_1, \dots, g_ig_{i+1}, \dots, g_{k+1}) & \hbox{for $0<i<k+1$}\\
\sigma(g_2,\dots,g_{k+1}) & \hbox{for $i=0$}
\end{array} \right.$$
\begin{definition} \label{defdoublecomplex}
For $\Gg$ a smooth \'{e}tale Leray groupoid, let's denote by
$\breve{C}^*(\Gg,\complex^\times(p)_\Gg)$ the total complex
\begin{eqnarray*}
\xymatrix{
\breve{C}^0(\Gg, \complex^\times(p)_\Gg)  \ar[r]^{\delta -d } &
         \breve{C}^1(\Gg,\complex^\times(p)_\Gg )
        \ar[r]^{\delta + d}  &  \breve{C}^2(\Gg, \complex^\times(p)_\Gg)
         \ar[r]^{\ \ \delta-d} & \cdots  \\
}
\end{eqnarray*}

induced
by the double complex
   $$\xymatrix{
     \vdots & \vdots & \vdots & & \vdots\\
     \breve{C}(\Gg_2, \complex^\times_\Gg) \ar[u]^\delta \ar[r]^{d \log} &
         \breve{C}(\Gg_2, \Aa^1_{\Gg,\complex})
        \ar[r]^d \ar[u]^\delta &  \breve{C}(\Gg_2, \Aa^2_{\Gg,\complex})
        \ar[u]^\delta \ar[r]^d & \cdots \ar[r]^d & \breve{C}(\Gg_2, \Aa^{p-1}_{\Gg,\complex})
         \ar[u]^\delta\\
      \breve{C}(\Gg_1, \complex^\times_\Gg) \ar[u]^\delta \ar[r]^{d \log} &
          \breve{C}(\Gg_1, \Aa^1_{\Gg,\complex})
        \ar[r]^d \ar[u]^\delta &  \breve{C}(\Gg_1, \Aa^2_{\Gg,\complex})
        \ar[u]^\delta \ar[r]^d & \cdots \ar[r]^d & \breve{C}(\Gg_1, \Aa^{p-1}_{\Gg,\complex})
        \ar[u]^\delta\\
      \breve{C}(\Gg_0, \complex^\times_\Gg) \ar[u]^\delta \ar[r]^{d \log} &
         \breve{C}(\Gg_0, \Aa^1_{\Gg,\complex})
          \ar[r]^d \ar[u]^\delta &  \breve{C}(\Gg_0, \Aa^2_{\Gg,\complex})
          \ar[u]^\delta \ar[r]^d & \cdots\ar[r]^d & \breve{C}(\Gg_0, \Aa^{p-1}_{\Gg,\complex})
          \ar[u]^\delta
    }$$
with $(\delta +(-1)^{i} d )$ as coboundary operator.
The {\v{C}}ech hypercohomology of the complex of sheaves
$\complex^\times(p)_\Gg$ is defined as the cohomology of the
{\v{C}}ech complex $\breve{C}(\Gg,\complex^\times(p)_\Gg)$:
$$\breve{H}^*(\Gg, \complex^\times(p)_\Gg):=H^*\breve{C}(\Gg,\complex^\times(p)_\Gg).$$
\end{definition}
As the $\Gg_i$'s are diffeomorphic to a disjoint union of
contractible sets -- Leray -- then the previous cohomology
actually matches the hypercohomology of the complex
$\complex^\times(p)_\Gg$, so we get
\begin{lemma} \label{cech=hyper}
The cohomology of the \u{C}ech complex
$\breve{C}^*(\Gg,\complex^\times(p)_\Gg)$ is isomorphic to the
hypercohomology of $\complex^\times(p)_\Gg$
$$\breve{H}^*(\Gg, \complex^\times(p)\Gg) \stackrel{\cong}{\to}
\hyper^*(\Gg,\complex^\times(p)_\Gg). $$
\end{lemma}

The Deligne cohomology groups classify the isomorphism classes of
what is known as $n$-gerbes with connective structure.
\begin{definition}
An $n$-gerbe with connective structure over $\Gg$ is an
$(n+1)$-cocycle of $\breve{C}^{n+1}(\Gg,\complex^\times(n+2)_\Gg)$.
Their isomorphism classes are classified by
$\hyper^{n+1}(\Gg,\complex^\times(n+2)_\Gg) =
\hyper^{n+2}(\Gg,\integer(n+2)^\infty_D)$.
\end{definition}

The following fact is more or less obvious, and relates this
definition with the one given by the authors in
\cite{LupercioUribe1}:

\begin{proposition}\label{olddefinition} To have a 1-gerbe  over $\Gg$ is the same thing
as to have a line bundle $\LL$  over $\Gg_1$ together with maps $\theta$, $h$
satisfying the
following properties:
\begin{itemize}
\item $\invers^* \LL \stackrel{\theta}{\cong} \LL^{-1}$
\item $\pi_1^* \LL \otimes \pi_2^* \LL \otimes \mult^*\invers^* \LL \stackrel{h}{\cong} 1$
\item $h: \Gg_1 \timests \Gg_1 \to \complex^\times$ is a 2-cocycle.
\end{itemize}
When the groupoid $\Gg$ is Leray, then the line bundle $\LL$ is trivial and all the information
is encoded in the 2-cocycle $h$. In this case a gerbe with connection will consist also of
a 1-form $A \in \Omega^1(\Gg_1)$, a 2-form $B \in
\Omega^2(\Gg_0)$ and a 3-form $K \in \Omega^3(\Gg_0)$ satisfying:
\begin{itemize}
\item $K=dB$
\item $\target^*B - \source^*B = dA$ and
\item $\pi_1^*A + \pi_2^*A -\mult^*A = -\sqrt{-1} h^{-1} dh$
\end{itemize}
\end{proposition}

As we will see via the holonomy map:
\begin{itemize}
\item A $0$-gerbe with connective structure induces a line bundle with connection over the groupoid $\Gg$ and a global 2-form on $\coarse{\Gg}$
\item A $1$-gerbe with connective structure induces what is known in the literature as a {\it gerbe with connection} over
$\Gg$ and a global 3-form on $\coarse{\Gg}$.
\end{itemize}

Before finishing this section let's point out that the group
$\hyper^{n-1}(\Gg,\complex^\times(n)_\Gg)$ is the only one that
encodes really new information as the following proposition clarifies.

\begin{proposition} \label{uninteresting part}
$$\hyper^{p}(\Gg,\integer(n)^\infty_D) \cong \hyper^{p-1}(\Gg,\complex^\times(n)_\Gg) = \left\{
\begin{array}{cc}
H^{p-1}(\Gg, \underline{\complex}^\times) = H^{p}(\Gg, \integer) & \mbox{for} \ p > n \\
H^{p-1}(\Gg, \complex^\times) & \mbox{for} \ p < n
\end{array} \right.$$
where $\underline{\complex}^\times$ stands for the sheaf of $\complex^\times$ valued functions. 
\end{proposition}

\begin{proof}
Let's have a look at the double complex of the definition
\ref{defdoublecomplex}. When $p>n$ the $p$-cocycles are over the
diagonal, so the information of all the columns of the double
complex besides the first is irrelevant. This is because the
sheaves $  \Aa^{i}_{\Gg,\complex}$ are acyclic. Now, when $p<n$
and $(h, \omega_1, \dots, \omega_{p-1})$ is a $(p-1)$-cocycle, by a
successive application of Poincar\'{e} lemma it is possible to find
an element $(f, \theta_1, \dots, \theta_{p-2})$ in
$\breve{C}^{p-2}(\Gg,\complex^\times(n)_\Gg)$ such that
$$(h, \omega_1, \dots, \omega_{p-1}) - (\delta +(-1)^{p-2}d)(f, \theta_1, \dots, \theta_{p-2}) = (h -\delta f, 0, \dots, 0)$$
with $d(h -\delta f)=0$, a locally constant $\complex^\times$ function. This implies the second isomorphism.
\end{proof}

After this brief summary of definitions we are ready to define the
holonomy map for {\bf{smooth \'{e}tale groupoids}}.

\section{Holonomy} \label{section Holonomy}

In the same way that a line bundle with connection over a manifold
$M$ induces a $\complex^\times$ valued function on the free loop
space of $M$, given by the holonomy around a loop, we can define
its analogous to smooth \'{e}tale groupoids. Let's recall that the
groupoid in mind is Leray, so we can make use of the {\v C}ech
description of the hypercohomology.

Let $W$ be and admissible cover of the circle associated to the set
$\{\alpha_0, \alpha_1, \dots, \alpha_n\}$ with $0=\alpha_0 < \alpha_1 < \cdots < \alpha_n=1$ as in section
\ref{sectionloopgroupoid}.

\begin{theorem} \label{holonomy}
There is a natural transgression map (holonomy)
$$\tau_1 :\breve{C}^1(\Gg,\complex^\times(2)_\Gg) \To \breve{C}^0(\Loop \Gg(W),\complex^\times_{\Loop \Gg(W)})$$
that sends cocycles to cocycles and that descends to cohomology
$$ \hyper^1\left(\Gg, \complex^\times_\Gg \stackrel{d \log}{\To} \Aa^1_{\Gg,\complex}\right)
 \To H^0\left(\Loop \Gg(W), \complex^\times_{\Loop \Gg(W)}\right).$$
\end{theorem}

\begin{proof}
First we will set up the notation. The pair $(h,A)$ will be an element in $\breve{C}^1(\Gg,\complex^\times(2)_\Gg)$
with $h : \Gg_1 \to \complex^\times$ and $A \in \Gamma(\Gg_0,\Aa^1_{\Gg,\complex})$. The object $\psi : \Ss^1_W \to \Gg$
of the loop groupoid $\Loop \Gg(W)$ will consist of maps $\psi_i : I_i=[\alpha_{i-1},\alpha_i] \to \Gg_0$ and arrows
$\psi :\{\alpha_1, \dots, \alpha_n\} \to \Gg_1$ such that
$$\source (\psi(\alpha_i)) = \psi_i(\alpha_i)  \ \ \ \  \ \ \ \ \  \mbox{and} \ \ \ \ \ \ \ \ \ \target (\psi(\alpha_i)) = \psi_{i+1}(\alpha_i)$$
and when $i=n$ $\target \psi(\alpha_n) = \psi_1(\alpha_0)$.

So $\tau_1(h,A)$ defines a function $H:\Loop \Gg(W)_0 \to \complex^\times$ as follows:
\begin{eqnarray}
H(\psi) := \exp\left( \sum_{i=1}^n \int_{I_i} \psi_i^* A \right) \prod_{i=1}^n h(\psi(\alpha_i))^{-1}. \label{definition of H}
\end{eqnarray}

It is clearly an homomorphism. We show now that $H$ descends to
cohomology. Suppose that $(h,A)$ is a 1-cocycle, its coboundary
$(d + \delta)(h,A)$ is zero, i.e.
\begin{eqnarray} \label{1-cocycle a}
\target^*A - \source^*A = - d \log h \ \ \ \mathrm{in} \ \ \ \Gg_1
\end{eqnarray}
\begin{eqnarray}
h(g_2) h(g_1g_2)^{-1}h(g_1) =1 \ \  \mbox{for}\ \ (g_1,g_2) \in \Gg_2, \label{1-cocycle b}
\end{eqnarray}
we want to see that the coboundary $\delta$ of $\tau_1(h,A)=H$ is also zero.
The cycle $\delta H$ is a function $\Loop \Gg(W)_1 \to \complex^\times$ that for the arrow
$\Lambda$ between $\psi$ and $\phi$ takes the value $\delta H (\Lambda) = H(\phi) H(\psi)^{-1}$.
The arrow $\Lambda$ will consist of maps $\Lambda_i : I_i \to \Gg_1$ such that
$$\Lambda_i(\alpha_i) \cdot \phi(\alpha_i) = \psi(\alpha_i) \cdot \Lambda_{i+1}(\alpha_i)$$
where $\Lambda_{n+1}(\alpha_n) := \Lambda_0(\alpha_0)$. In the following diagram the dark lines are the images of the intervals
in $\Gg_0$ and the arrows are elements in $\Gg_1$.
\begin{eqnarray}
\includegraphics[height=3in]{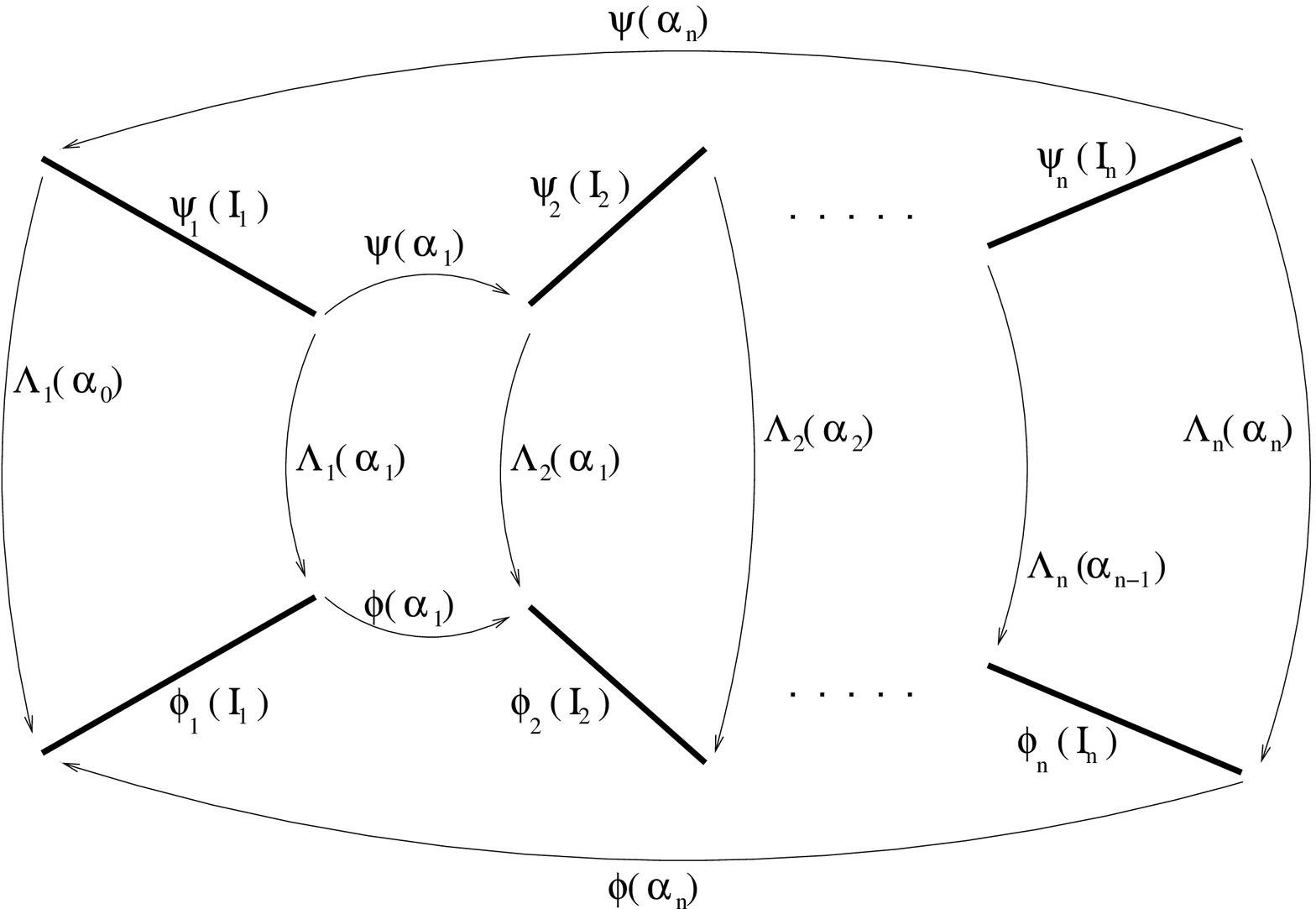} \label{graph morphism}
\end{eqnarray}

Making use of  the property \ref{1-cocycle a}
of the 1-cocycle $(h,A)$ we get the following set of equalities:
\begin{eqnarray*}
\frac{\exp\left( \int_{I_i}  \phi_i^*A \right)}
{\exp\left( \int_{I_i} \psi_i^*A \right)}
& = &  \exp\left( \int_{I_i} \Lambda_i^*( \target^* A - \source^* A) \right)  \\
& = & \exp\left( \int_{I_i}  \Lambda_i^*(-d \log h) \right)  \\
& =  & \frac{h(\Lambda_i(\alpha_{i-1}))}{h(\Lambda_i(\alpha_i))},
\end{eqnarray*}
and using property \ref{1-cocycle b} we have
\begin{eqnarray*}
\delta H (\Lambda) = \frac{H(\phi)}{H(\psi)} & = & \prod_{i=1}^n \frac{h(\Lambda_i(\alpha_{i-1}))}{h(\Lambda_i(\alpha_i))}
\frac{h(\psi(\alpha_i))}{h(\phi(\alpha_i))} \\
& = & \prod_{i=1}^n \frac{ h(\psi(\alpha_i)) h( \Lambda_{i+1}(\alpha_i))}{h(\Lambda_i(\alpha_i)) h (\phi(\alpha_i))}\\
& = & \prod_{i=1}^n \frac{h \left(\Lambda_i(\alpha_i) \cdot \phi(\alpha_i) \right)}{h \left( \psi(\alpha_i) \cdot \Lambda_{i+1}(\alpha_i)\right)} =1
\end{eqnarray*}
This means that $H$ is invariant under the action of $\Loop \Gg(W)_1$ therefore it defines  a map
$$H : \Loop \Gg(W) / \sim \to \complex^\times$$

Now if $(h,A)$ is a coboundary, i.e. $(h,A)= (\delta f, -d \log f)$ for some $f : \Gg_0 \to \complex^\times$ then
$H = \tau_1 (\delta f, -d \log f)$ will become
\begin{eqnarray*}
H(\psi) & = & \exp\left( \sum_{i=1}^n \int_{I_i} \psi_i^* (-d \log f) \right) \prod_{i=1}^n \delta f(\psi(\alpha_i))^{-1} \\
& = & \prod_{i=1}^n \frac{f(\psi_i(\alpha_{i-1}))}{f(\psi_i(\alpha_{i}))}
\prod_{n=1}^n \frac{f(\source(\psi(\alpha_i)))}{f(\target(\psi(\alpha_i)))} =1.
\end{eqnarray*}

Hence the map $\tau_1$ descends to cohomology.

\end{proof}

Another way to understand the previous result is the following.
The pair $(h,A)$ represents a complex line bundle  with connection
$(L, \Delta)$ over the groupoid $\Gg$. The function $\tau_1(h,A)$
assigns a complex number to every element $\psi$ in the loop
groupoid. This number represents an endomorphism of the fiber
$L_\psi$ obtained through the parallel transport given by the
connection $\Delta$.

If we now take a refinement $W'$ of the cover $W$ associated to
the set \linebreak $\{\alpha_0, \dots, \alpha_{i-1}, \beta,
\alpha_{i}, \dots, \alpha_n\}$, $\rho^{W'}_W : \Ss^1_{W'} \to
\Ss^1_W$ the natural morphism and $\psi : \Ss^1_W \to \Gg$ a loop,
we can see that for $\psi' := \psi \circ \rho^{W'}_W$ the equality
$H(\psi) = H(\psi')$ holds. This because the morphism
$\psi(\beta)$ is equal to $\ident \psi_i(\beta)$ and by property
\ref{1-cocycle b} the function $h$ restricted to $\ident(\Gg_0)$
is constant and equal to $1$. Then the function $H$ is stable
under cover refinement, and therefore we can take the inverse
limits on the admissible covers and to obtain the holonomy morphism
over the loop groupoid $\Loop \Gg$:

\begin{proposition}
There is a natural transgression map (holonomy)
$$\tau_1 :\breve{C}^1(\Gg,\complex^\times(2)_\Gg) \To \breve{C}^0(\Loop \Gg,\complex^\times_{\Loop \Gg})$$
that sends cocycles to cocycles and that descends to cohomology
$$ \hyper^1\left(\Gg, \complex^\times_\Gg \stackrel{d \log}{\To} \Aa^1_{\Gg,\complex}\right)
 \To H^0\left(\Loop \Gg, \complex^\times_{\Loop \Gg}\right).$$
\end{proposition}

\begin{rem}
Using the definition of the holonomy given by Brylinski
\cite{Brylinski} in Lemma 6.1.2 and taking $\Gg$ to be  a Leray
groupoid naturally associated to a manifold $M$  as in Example
\ref{manifoldgroupoid} (i.e. $\Gg$ is built out of a open
contractible cover of $M$) we see that the previous map matches
the holonomy of a connection in a line bundle around a loop in $M$.
\end{rem}

\section{The Line bundle over the loop groupoid} \label{section Line bundle over the loop groupoid}

From a gerbe with connection over the groupoid $\Gg$  we are going
to construct a line bundle over the loop groupoid, in a way that
is compatible with  the transgression map on a manifold. The main
result of this section is:

\begin{theorem} \label{theoremlinebundle}
There is a natural homomorphism
$$\tau_2 :\breve{C}^2(\Gg,\complex^\times(3)_\Gg) \To \breve{C}^1(\Loop \Gg(W),\complex^\times(2)_{\Loop \Gg(W)})$$
that sends 2-cocycles to 1-cocycles (i.e. gerbes with connection
over $\Gg$ to line bundles with connection over the loop groupoid
$\Loop \Gg$),
  commutes with the coboundary operator (i.e. $\tau_2 \circ (\delta + d) =(\delta - d ) \circ \tau_1$)
and therefore induces a map in cohomology
$$\hyper^2\left(\Gg, \complex_\Gg^\times \stackrel{d \log}{\to} \Aa^1_{\Gg,\complex} \stackrel{d}{\to} \Aa^2_{\Gg, \complex} \right)
\To \hyper^1\left(\Loop \Gg(W), \complex_{\Loop \Gg(W)}^\times
\stackrel{d\log}{\to} \Aa^1_{\Loop \Gg(W), \complex} \right).$$
\end{theorem}

\begin{proof}
Let's first fix the notation. The triple $(h,A,B)$ will be an element of $ \breve{C}^2(\Gg,\complex^\times(3)_\Gg)$ with
$h: \Gg_2 \to \complex^\times$, $A \in \Gamma(\Gg_1,\Aa^1_{\Gg,\complex})$  and $B \in \Gamma(\Gg_0,\Aa^2_{\Gg,\complex})$. The
arrow $\Lambda$ of the loop groupoid $\Loop \Gg(W)$ between the objects $\psi$ and $\phi$ will be defined as in theorem
\ref{holonomy}. The arrow $\nu$ of the tangent loop groupoid $T_\Lambda \Loop \Gg(W)$ between the objects $\xi \in T_\psi \Loop \Gg(W)$
 and $\zeta\in T_\phi \Loop \Gg(W)$ (as in definition \ref{tangent loop groupoid}) will consist of maps
$\xi_i,\zeta_i : I_i \to T\Gg_0$, $\nu_i : I_i \to T\Gg_1$ and $\xi,\zeta : \{\alpha_1, \dots, \alpha_n\} \to T\Gg_1$
such that:
$$\source (\xi_i(\alpha_i)) = \xi_i(\alpha_i), \ \ \target (\xi(\alpha_i)) = \xi_{i+1}(\alpha_i), \ \
\source (\zeta(\alpha_i)) = \zeta_i(\alpha_i), \ \ \target (\zeta(\alpha_i)) = \zeta_{i+1}(\alpha_i)$$
$$\nu_i(\alpha_i) \cdot \zeta(\alpha_i) = \xi(\alpha_i) \cdot \nu_{i+1}(\alpha_i)$$

and when $i=n$ $\target \xi(\alpha_n) = \xi_1(\alpha_0)$, $\target \zeta(\alpha_n) = \zeta_1(\alpha_0)$
 and $\nu_{n+1}(\alpha_n) := \nu_0(\alpha_0)$.

Now we are ready to define $\tau_2(h,A,B)$. It will consist of the pair $(F,\Delta)
\in \breve{C}^1(\Loop \Gg(W),\complex^\times(1)_{\Loop \Gg(W)})$ with $F : \Loop \Gg (W) \to \complex^\times$ a map,
 and $\Delta : T \Loop \Gg(W)_0 \to \complex$ a linear functional on the tangent loop space.

\begin{definition}
\begin{eqnarray*} \label{definitionF}
F(\Lambda):= \exp \left( \sum_{i=1}^n \int_{I_i} \Lambda_i^* A \right) \prod_{i=1}^n
\frac{h(\psi(\alpha_i),\Lambda_{i+1}(\alpha_i))}{h(\Lambda_i(\alpha_i),\phi(\alpha_i))}
\end{eqnarray*}
\begin{eqnarray*} \label{definitionDelta}
\langle \Delta_\psi , \xi \rangle := \sum_{i=1}^n \int_{I_i} B \left( \frac{d \psi_i}{dt}, \xi_i(t) \right) dt
+ \sum_{i=1}^n \langle A_{\psi(\alpha_i)} , \xi(\alpha_i) \rangle
\end{eqnarray*}
\end{definition}

In this way the map $\tau_2$ is clearly an homomorphism.

The other two statements of the theorem will be proven separately.

\begin{proposition}
$\tau_2$ sends cocycles to cocycles.
\end{proposition}

\begin{proof}
In other words we need to prove that $(\delta -d)(h,A,B)=0$ implies $(\delta +d)(F,\Delta) =0$, i.e.
a gerbe with connection over the groupoid induces a line bundle with connection over the loop groupoid.

The cocycle condition $(\delta h, \delta A - d \log h, \delta B - dA )=0$ implies:
\begin{eqnarray}
h(a,b) h(a,bc)^{-1} h(ab,c)h(b,c)^{-1} & =& 1 \ \ \ \mathrm{for} \  \ \ \ (a,b,c) \in \Gg_3 \label{2-cocycle h} \\
\pi_2^*A + \pi_1^*A - \mult^*A & =& d \log h \ \ \ \ \mathrm{in}  \ \ \ \ \Gg_2 \label{2-cocycle A} \\
\target^* B - \source^* B &=& d A \ \ \ \ \mathrm{in} \ \ \ \ \Gg_1. \label{2-cocycle B}
\end{eqnarray}

Let's prove first that $\delta F = 1$. This in particular implies that the map $F : \Loop \Gg (W) \to \complex^\times$
is a morphism of groupoids.

Let $\Lambda$ and $\Omega$ be two arrows in the loop groupoid with $\psi \stackrel{\Lambda}{\to} \phi \stackrel{\Omega}{\to} \gamma$,
we need to calculate $\delta F(\Lambda,\Omega) = F(\Lambda)F(\Omega)F(\Lambda \cdot \Omega)^{-1}$.

Using the property \ref{2-cocycle A} we have that
$$\exp\left(\int_{I_i} \Lambda_i^* A + \Omega_i^* A - (\Lambda \cdot \Omega)_i^* A \right)=
\frac{h(\Lambda_i(\alpha_i),\Omega_i(\alpha_i))}{ h(\Lambda_i(\alpha_{i-1}),\Omega_i(\alpha_{i-1}))} $$
and by applying property \ref{2-cocycle h} to the triples
$$(\psi(\alpha_i), \Lambda_{i+1}(\alpha_i), \Omega_{i+1}(\alpha_i)) \ \ \ \mbox{and} \ \ \
(\Lambda_i(\alpha_i),\Omega_i(\alpha_i), \gamma(\alpha_i))$$
\begin{eqnarray}
\includegraphics[height=2.6in]{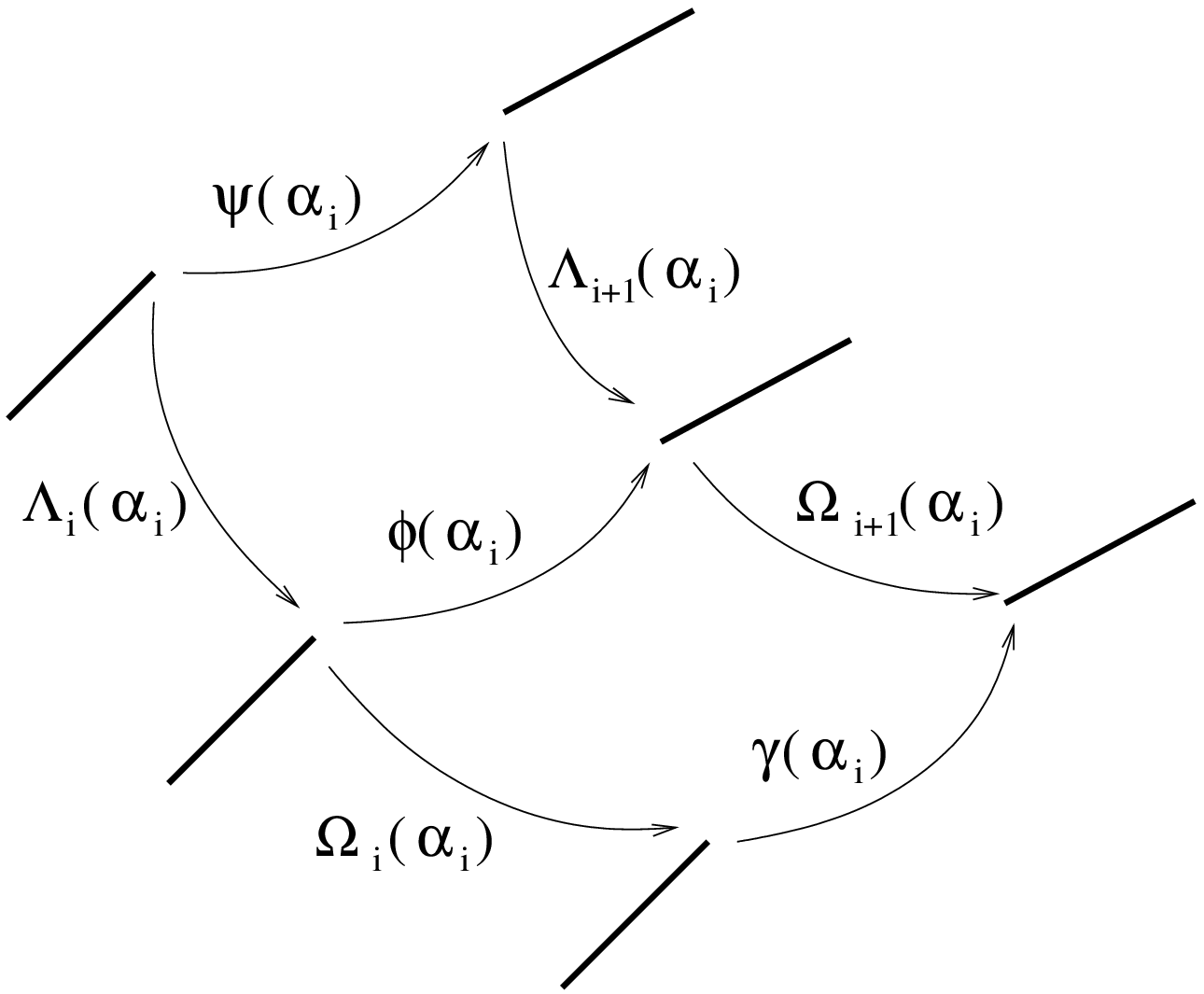} \label{graph composition}
\end{eqnarray}

we get
\begin{eqnarray*}
\frac{h(\psi(\alpha_i),\Lambda_{i+1}(\alpha_i))}{h(\Lambda_i(\alpha_i),\phi(\alpha_i))}
\frac{h(\phi(\alpha_i),\Omega_{i+1}(\alpha_i))}{h(\Omega_i(\alpha_i),\gamma(\alpha_i))}
\left(\frac{h(\psi(\alpha_i),(\Lambda\cdot \Omega)_{i+1}(\alpha_i))}{h((\Lambda\cdot \Omega)_i(\alpha_i),\gamma(\alpha_i))} \right)^{-1}
=  \\ \frac{h(\Lambda_{i+1}(\alpha_i),\Omega_{i+1}(\alpha_i))}{h(\Lambda_i(\alpha_i),\Omega_i(\alpha_i))}
\end{eqnarray*}
Multiplying the last two equations and making the product over the $i$'s it follows that $\delta F(\Lambda, \Omega)=1$.

Now we will prove that $(\delta \Delta + d \log F) =0$. Let $\nu \in T_\Lambda \Loop \Gg(W)$, we want to check
$$\langle - (d \log F)_\Lambda , \nu \rangle = \langle (\delta \Delta)_\Lambda , \nu \rangle$$
and we will do so by integrating over a 1-parameter thickening of
the path $\Lambda$ in he direction of $\nu$. Via the tubular
neighborhood  diffeomorphism and the fact that $\nu$ and $\Lambda$
are determined by a finite number of maps over compact sets, we
can find a one-parameter family $\Lambda^s \in \Loop \Gg(W)$, with
$s \in [-\epsilon, \epsilon]$ for $\epsilon$ sufficiently small,
such that $\Lambda = \Lambda^0$ and $\frac{d \Lambda^s}{ds} =
\nu$. We claim that
$$\int_{-\epsilon}^\epsilon \langle - (d \log F)_{\Lambda^s} , \frac{d \Lambda^s}{ds} \rangle ds =
\int_{-\epsilon}^\epsilon \langle (\delta \Delta)_{\Lambda^s}, \frac{d \Lambda^s}{ds} \rangle ds $$

Let's first elaborate on the left hand side (LHS). The steps will be outlined after the set of equalities.

\begin{eqnarray}
LHS & = & - \sum_{i=1}^n \left( \int_{I_i} \left( \Lambda^\epsilon_i \right)^*A - \left( \Lambda^{-\epsilon}_i \right)^*A \right)
\label{LHS1} \\
& & - \int_{-\epsilon}^\epsilon \left\langle d \log \prod_{i=1}^n
\frac{h(\psi^s(\alpha_i),\Lambda^s_{i+1}(\alpha_i))}{h(\Lambda_i^s(\alpha_i),\phi^s(\alpha_i))}, \frac{d \Lambda^s}{ds}
\right\rangle ds \label{LHS2} \\
& = &  - \sum_{i=1}^n \left( \int_{I_i} \langle A_{\Lambda^\epsilon_i(t)}, \frac{d \Lambda^\epsilon_i}{dt} \rangle -
\langle A_{\Lambda^{-\epsilon}_i(t)}, \frac{d \Lambda^{-\epsilon}_i}{dt} \rangle dt \right) \label{LHS3}\\
& & - \sum_{i=1}^n  \int_{-\epsilon}^\epsilon \left( \langle A_{\psi^s(\alpha_i)}, \frac{d \psi^s}{ds}(\alpha_i) \rangle -
\langle A_{\phi^s(\alpha_i)}, \frac{d \phi^s}{ds}(\alpha_i) \rangle \right. \label{LHS4}\\
& & \left. + \langle A_{\Lambda^s_{i+1}(\alpha_i)}, \frac{d \Lambda^s_{i+1}}{ds}(\alpha_i) \rangle -
\langle A_{\Lambda^s_{i}(\alpha_i)}, \frac{d \Lambda^s_{i}}{ds}(\alpha_i) \rangle \right) ds \label{LHS5}
\end{eqnarray}
Line \ref{LHS1} is obtained after evaluating the integral at the end points $-\epsilon$ and $\epsilon$.
Line \ref{LHS3} is the same as line \ref{LHS1} but written in a different way, and lines \ref{LHS4} and \ref{LHS5}
are obtained from line \ref{LHS2} after using property \ref{2-cocycle A} and the fact that
$$\psi^s(\alpha_i) \cdot \Lambda_{i+1}^s(\alpha_i) = \Lambda_i^s(\alpha_i) \cdot \phi^s(\alpha_i).$$

For the right hand side we need to make use of Stoke's theorem.

\begin{eqnarray}
RHS & = & \int_{-\epsilon}^\epsilon \langle \Delta_{\phi_i^s}, \frac{d \phi_i^s}{ds} \rangle -
 \langle \Delta_{\psi_i^s}, \frac{d \psi_i^s}{ds} \rangle ds \label{RHS1} \\
& = &  \int_{-\epsilon}^\epsilon \left( \sum_{i=1}^n \int_{I_i} B\left( \frac{d \phi_i^s}{dt} , \frac{d \phi_i^s}{ds} \right)
- B\left( \frac{d \psi_i^s}{dt} , \frac{d \psi_i^s}{ds} \right) dt \right) ds \label{RHS2} \\
& &  +  \sum_{i=1}^n  \int_{-\epsilon}^\epsilon \left( \langle A_{\phi^s(\alpha_i)}, \frac{d \phi^s}{ds}(\alpha_i) \rangle -
\langle A_{\psi^s(\alpha_i)}, \frac{d \psi^s}{ds}(\alpha_i) \rangle \right) ds \label{RHS3} \\
& = & \sum_{i=1}^n \int_{-\epsilon}^\epsilon \int_{I_i} dA\left( \frac{d \Lambda_i^s}{dt} , \frac{d \Lambda_i^s}{ds} \right) dt \ ds
\label{RHS4} \\
& &  +  \sum_{i=1}^n  \int_{-\epsilon}^\epsilon \left( \langle A_{\phi^s(\alpha_i)}, \frac{d \phi^s}{ds}(\alpha_i) \rangle -
\langle A_{\psi^s(\alpha_i)}, \frac{d \psi^s}{ds}(\alpha_i) \rangle \right) ds \label{RHS5} \\
& = & \sum_{i=1}^n \left( \int_{I_i} \langle A_{\Lambda^{-\epsilon}_i(t)}, \frac{d \Lambda^{-\epsilon}_i}{dt} \rangle -
\langle A_{\Lambda^{\epsilon}_i(t)}, \frac{d \Lambda^{\epsilon}_i}{dt} \rangle dt \right) \label{RHS6}\\
& & + \sum_{i=1}^n  \int_{-\epsilon}^\epsilon \left(
  \langle A_{\Lambda^s_{i}(\alpha_i)}, \frac{d \Lambda^s_{i}}{ds}(\alpha_i) \rangle -
\langle A_{\Lambda^s_{i+1}(\alpha_i)}, \frac{d \Lambda^s_{i+1}}{ds}(\alpha_i) \rangle \right) ds \label{RHS7} \\
& &  +  \sum_{i=1}^n  \int_{-\epsilon}^\epsilon \left( \langle A_{\phi^s(\alpha_i)}, \frac{d \phi^s}{ds}(\alpha_i) \rangle -
\langle A_{\psi^s(\alpha_i)}, \frac{d \psi^s}{ds}(\alpha_i) \rangle \right) ds \label{RHS8}
\end{eqnarray}
Line \ref{RHS1} is obtained after applying the coboundary operator to $\Delta$. Expanding via the definition of $\Delta$
we get lines \ref{RHS2} and \ref{RHS3}. From property \ref{2-cocycle B} we get line \ref{RHS4} from \ref{RHS2}. And lines
\ref{RHS5} and \ref{RHS6} come from line \ref{RHS4} and Stoke's theorem (evaluating $A$ at the boundary).

Lines \ref{LHS3}, \ref{LHS4} and \ref{LHS5}  match lines \ref{RHS6}, \ref{RHS8} and \ref{RHS7} respectively. Therefore
LHS = RHS and the pair $(F, \Delta)$ is a 1-cocycle.
\end{proof}

The second part of the theorem states

\begin{proposition}
The map $\tau_2$
commutes with the coboundary operator (i.e. $\tau_2 \circ (\delta + d) =(\delta - d ) \circ \tau_1$)
\end{proposition}
\begin{proof}

Let the pair $(f,G)$ be in $\breve{C}^1(\Gg,\complex^\times(2)_\Gg)$, $H := \tau_1(f,G)$, $(\delta - d)H = (\delta H, -d \log H)$,
$(\delta + d)(f,G) =
(\delta f, \delta G + d \log f , dG)$ and $(F,\Delta) := \tau_2(\delta f, \delta G + d \log f , dG)$.
We want to prove that $$(\delta H, -d \log H) = (F,\Delta).$$

Replacing in $F(\Lambda)$ we get:
\begin{eqnarray*}
F(\Lambda) & = & \exp \left( \sum_{i=1}^n \int_{I_i} \Lambda_i^*(\delta G + d \log f) \right)
\prod_{i=1}^n \frac{\delta f (\psi(\alpha_i), \Lambda_{i+1}(\alpha_i))}{\delta f(\Lambda_i(\alpha_i), \phi(\alpha_i))} \\
&=& \exp   \left( \sum_{i=1}^n \int_{I_i} \phi^* G - \psi^*G \right) \prod_{i=1}^n
\frac{f(\Lambda_i(\alpha_i))}{f(\Lambda_i(\alpha_{i-1}))}
\prod_{i=1}^n \frac{ f (\psi(\alpha_i))f (\Lambda_{i+1}(\alpha_i))}{ f(\Lambda_i(\alpha_i)) f(\phi(\alpha_i))}\\
& = &  \exp   \left( \sum_{i=1}^n \int_{I_i} \phi^* G - \psi^*G  \right)
 \frac{ f (\psi(\alpha_i)))}{  f(\phi(\alpha_i))} \\
&=& H(\phi) H(\psi)^{-1} = \delta H (\Lambda)
\end{eqnarray*}

To prove $-d \log H = \Delta$ we proceed as in the previous proposition. We claim that
$$\int_{-\epsilon}^\epsilon \langle - (d \log H)_{\psi^s} , \frac{d \psi^s}{ds} \rangle ds =
\int_{-\epsilon}^\epsilon \langle ( \Delta)_{\psi^s}, \frac{d \psi^s}{ds} \rangle ds. $$

Elaborating on each side in the same way it was done for the previous proposition we have that

\begin{eqnarray}
LHS & = & - \sum_{i=1}^n \left( \int_{I_i} ( \psi^\epsilon_i)^* G - ( \psi^{-\epsilon}_i)^* G \right) \label{lhs1}\\
& & - \sum_{i=1}^n \int_{-\epsilon}^\epsilon \langle (d \log f^{-1})_{\psi^s(\alpha_i)} , \frac{d \psi^s}{ds}(\alpha_i) \rangle ds \label{lhs2}
\end{eqnarray}
\begin{eqnarray}
RHS & = &\sum_{i=1}^n  \int_{-\epsilon}^\epsilon \int_{I_i} dG \left(\frac{d \psi^s}{dt},\frac{d \psi^s}{ds} \right) dt \ ds
\label{rhs1} \\
& & \sum_{i=1} \int_{-\epsilon}^\epsilon \langle (\delta G + d \log f)_{\psi^s(\alpha_i)} , \frac{d \psi^s}{ds}(\alpha_i) \rangle ds
\label{rhs2} \\
& = & \sum_{i=1}^n \left( \int_{I_i} ( \psi^{-\epsilon}_i)^* G - ( \psi^{\epsilon}_i)^* G \right) \label{rhs3}\\
& & +\sum_{i=1}^n \int_{-\epsilon}^\epsilon \langle G_{\psi_i^s(\alpha_i)}, \frac{d \psi_i^s}{ds}(\alpha_i) \rangle -
\langle G_{\psi_i^s(\alpha_{i-1})}, \frac{d \psi_i^s}{ds}(\alpha_{i-1}) \rangle ds \label{rhs4}\\
& & +\sum_{i=1}^n \int_{-\epsilon}^\epsilon \langle G_{\psi_{i+1}^s(\alpha_i)}, \frac{d \psi_{i+1}^s}{ds}(\alpha_i) \rangle -
\langle G_{\psi_i^s(\alpha_{i})}, \frac{d \psi_i^s}{ds}(\alpha_{i}) \rangle ds \label{rhs5}\\
& & + \sum_{i=1}^n \int_{-\epsilon}^\epsilon \langle (d \log f)_{\psi^s(\alpha_i)} , \frac{d \psi^s}{ds}(\alpha_i) \rangle ds \label{rhs6}
\end{eqnarray}

Where the lines \ref{lhs1}, \ref{lhs2}, \ref{rhs1} and \ref{rhs2} are obtained after replacing the given information and lines
\ref{rhs3} and \ref{rhs4} come from \ref{rhs1} and Stoke's theorem. Lines \ref{rhs4} and \ref{rhs5}  are equal 
with opposite signs; therefore we have
that $LHS = RHS.$

So we have that the following square is commutative:

\begin{eqnarray*}
\xymatrix{
\breve{C}^1(\Gg,\complex^\times(2)_\Gg) \ar[r]^{\tau_1} \ar[d]_{(\delta +d)} &
 \breve{C}^0(\Loop \Gg(W),\complex^\times_{\Loop \Gg(W)}) \ar[d]^{(\delta -d)} \\
\breve{C}^2(\Gg,\complex^\times(3)_\Gg) \ar[r]_{\tau_2} & \breve{C}^1(\Loop \Gg(W),\complex^\times(2)_{\Loop \Gg(W)}) \\
}
\end{eqnarray*}

\end{proof}
With the previous two propositions it is clear that $\tau_2$ descends to a map in cohomology and theorem \ref{theoremlinebundle}
follows.

\end{proof}

Taking the limit over the admissible covers we obtain the following statement
\begin{cor}There is a natural homomorphism
$$\tau_2 :\breve{C}^2(\Gg,\complex^\times(3)_\Gg) \To \breve{C}^1(\Loop \Gg,\complex^\times(2)_{\Loop \Gg})$$
that sends 2-cocycles to 1-cocycles (i.e. gerbes with connection over $\Gg$ to line bundles with connection over the loop groupoid),
  commutes with the coboundary operator (i.e. $\tau_2 \circ (\delta + d) =(\delta - d ) \circ \tau_1$)
and therefore induces a map in cohomology
$$\hyper^2\left(\Gg, \complex_\Gg^\times \stackrel{d \log}{\to} \Aa^1_{\Gg,\complex} \stackrel{d}{\to} \Aa^2_{\Gg, \complex} \right)
\To \hyper^1\left(\Loop \Gg, \complex_{\Loop \Gg}^\times
\stackrel{d\log}{\to} \Aa^1_{\Loop \Gg, \complex} \right).$$
\end{cor}

\subsection{Manifolds} \label{subsection Manifolds}
In the case that the groupoid $\Mm$ represents a manifold $M$ via
a good open cover $\{U_\alpha\}_{\alpha \in I}$ (by good we mean
that the finite intersections $U_{\alpha_1 \dots
\alpha_n}:=U_{\alpha_1} \cap \cdots \cap U_{\alpha_n}$ are either
empty or contractible -namely a Leray cover-), with
$$\Mm_0 = \bigsqcup_{\alpha \in I} U_\alpha \ \ \ \ \mbox{   and   } \ \ \ \  \Mm_1 = \bigsqcup_{(\alpha,\beta) \in I^2} U_{\alpha\beta}$$
 we obtain the construction introduced by
Gaw{\c e}dzki  of a line bundle with connective structure over
$\LL M$ via a gerbe with connection over $M$ \cite[Pag.
108-113]{Gawedzki} (it can also be found in Brylinski's book
\cite[Prop. 6.5.1]{Brylinski}). The information of the cocycle $(h,A,B)$ whose cohomology class lies on
 $\hyper^2\left(\Mm, \complex_\Mm^\times \stackrel{d \log}{\to} \Aa^1_{\Mm,\complex} \stackrel{d}{\to} \Aa^2_{\Mm, \complex} \right)$
is equivalent to the data
$$h_{\alpha\beta\gamma} : U_{\alpha\beta\gamma} \to \complex^\times, \ \ \ \  A_{\alpha\beta} \in \Omega^2(U_{\alpha\beta})\otimes \complex
\ \ \ \mbox{and} \ \  \ \ B_\alpha \in \Omega^1(U_\alpha) \otimes
\complex$$ so that
\begin{eqnarray*}
h_{\alpha\beta\gamma}h^{-1}_{\alpha\beta\delta}h_{\alpha\gamma\delta}h^{-1}_{\beta\gamma\delta} & =&1  \ \ \ \mbox{in} \ \ \ U_{\alpha\beta\gamma\delta}\\
A_{\alpha\beta} +A_{\beta\gamma}-A_{\alpha\gamma} &=& d \log h_{\alpha\beta\gamma} \ \ \ \ \mbox{in} \ \ \  U_{\alpha\beta\gamma}\\
B_\beta - B_\alpha &=& d A_{\alpha\beta} \ \ \ \  \mbox{in} \ \ \
U_{\alpha\beta}.
\end{eqnarray*}
Let $\Lambda : \psi \to \phi$ be an arrow in the loop groupoid
$\Loop \Mm (W)$.
Recall that as $M$ is a manifold the category $\Mm$ does not have
automorphisms besides the identity (i.e. there is only one arrow from a point to itself) .
Define the indices $\kappa_i, \lambda_i \in I$ such that
$$\psi_i (I_i) \subset U_{\kappa_i} \ \ \ \ \phi_i (I_i) \subset U_{\lambda_i},$$
and therefore $\Lambda(I_i) \subset U_{\kappa_i \lambda_i}$.
Hence, the formula \ref{definitionF} can be written as:
$$F(\Lambda) = \exp \left( \sum_{i=1}^n \int_{I_i} \Lambda_i^*A_{\kappa_i\lambda_i} \right)
\prod_{i=1}^n
\frac{h_{\kappa_i\kappa_{i+1}\lambda_{i+1}}(\Lambda(\alpha_i))}{h_{\kappa_i\lambda_{i}\lambda_{i+1}}(\Lambda(\alpha_i))}$$
with $\kappa_{n+1} =\kappa_1$ and $\lambda_{n+1} =\lambda_1$.
 If $\xi \in T_\psi \Loop \Mm$ is a vector field over $\psi$
(a tangent vector of the $\Loop \Mm$), then the formula
\ref{definitionDelta} can be written as:
$$(\Delta, \xi)_{\psi} =  \sum_{i=1}^{n} \int_{I_i} B_{\kappa_i}\left( \frac{d \psi_i}{dt}, \xi_i(t)\right)dt +
\sum_{i=1}^{n} \langle A_{\kappa_i\kappa_{i+1}}(\psi(\alpha_i)),
\xi(\alpha_i) \rangle. $$ This assignment matches the ones given by
Gaw{\c e}dzki \cite[Pag. 111]{Gawedzki} and Brylinski \cite[Pag.
250]{Brylinski}. As $\LL M $, the loop space of
$M$, and the loop groupoid $\Loop \Mm$ are Morita equivalent (see
\cite[Prop. 5.1.3]{LupercioUribe2}) we can deduce  Brylinski's
result:
\begin{proposition} \cite[Prop. 6.5.1]{Brylinski}
The assignment $(h,A,B) \mapsto (F,\Delta)$ induces a group
homomorphisms
$$\hyper^2\left(M, \complex_M^\times \stackrel{d \log}{\to} \Aa^1_{M,\complex} \stackrel{d}{\to} \Aa^2_{M, \complex} \right)
\To \hyper^1\left(\LL M, \complex_{\LL M}^\times
\stackrel{d\log}{\to} \Aa^1_{\LL M, \complex} \right)$$ and is
equal to the opposite of the transgression map from sheaves of
groupoids with connective structure and curving over $M$ to line
bundles with connection over $\LL M$.
\end{proposition}

\subsection{Global quotients} \label{subsection Global quotients}

For the purpose of illustration let us consider an orbifold of the
form $[M/G]$ obtained from a manifold in which $G$ acts as a
finite subgroup of $\Diff(M)$. Moreover let us assume that the
gerbe $(h,A,B)$  that we consider can be represented on the groupoid
$\Xx$ with morphisms $\Xx_1 = M \times G$ and objects $\Xx_0=M$,
where the arrow $(m,g)$ takes the object $m$ to the object $mg$.
Namely:
\begin{itemize}
   \item $B \in \Omega^2(M) \otimes \complex$
   \item $A \in \Omega^1(M \times G) \otimes \complex $
   \item $h \colon \Xx_2 = M \times G \times G \rightarrow
   \complex^\times$,
\end{itemize}
and writing $A_g := A|_{M \times \{g\}}$ and $h_{g,k} := h|_{M \times \{g\} \times \{k\}}$ for $g,k \in G$ ,we have
\begin{eqnarray*}
g^* B - B &=& d A_g \\
A_g +A_k + A_{(gk)^{-1}} &=& d \log h_{g,h}.
\end{eqnarray*}

Let $\Lambda =(\phi,g,k)$ be an arrow of the loop groupoid $\Loop \Xx$ (as explained in example \ref{example loop groupoid})
 from $(\phi,g)$ to $(\phi \cdot k, k^{-1}gk)$ where $\phi : [0,1] \to M$ and $\phi(0)g = \phi(1)$. Then
 the formula \ref{definitionF} can be written in this
case as:

$$F(\Lambda)= \exp \left( \int_0^1 \phi^* A_k \right) \frac{h_{g,k}(\phi(0))}{h_{k,k^{-1}gk}(\phi(0))}.$$

And if $\xi$ is a vector field along $\phi$ (i.e. $\xi \in \Gamma([0,1], \phi^*TM)$) with $\xi(0)g = \xi(1)$
then the functional of equation \ref{definitionDelta} can be expressed as:

$$(\Delta, \xi)_\phi =  \int_0^1 B \left( \frac{d \phi}{dt},  \xi(t) \right)dt  + \langle A_g(\phi(0)), \xi(0) \rangle$$

\subsection{Discrete torsion} \label{subsection Discrete torsion}

When considering conformal field theories on an orbifold $[M/G]$ it
is well-known that for any non-trivial cocycle $ \varepsilon : G
\times G \to \complex^\times$ and $[\alpha] \in H^2(G,
\complex^\times)$,
 a new model can be defined
by weighting the twisted sectors of the orbifold with a
non-trivial phase, the so-called discrete torsion (see
\cite{Vafa}). As argued in \cite{LupercioUribe3} this can be
seen as a choice of flat $B$-field over the target stack $[M/G]$.

Let $\Xx$ be a groupoid associated to $[M/G]$, and as in the previous section
take $\Xx_0:=M$ and $\Xx_1:=M\times G$ with the natural source and
target maps. Let $\bar{G} := * \times G \twoarrows *$ be the
natural groupoid representative of $G$. The morphism $\Xx \to
\bar{G}$ induces a monomorphism $H^2(G, \complex^\times) \to
\hyper^2 (\Xx, \complex^\times(3)_\Xx)$ that allows to define a
flat gerbe as follows
\begin{eqnarray*}
h: \Xx_2 = M \times G \times G & \to & \complex^\times \\
(x,g_1, g_2) & \mapsto & \varepsilon(g_1,g_2)\\
A=B=0.
\end{eqnarray*}

As the flat gerbe only depends on the group $G$ and not on the
geometry of $M$, we need not to work with an open cover of the
orbifold $\Xx$. We know \cite[Prop. 6.1.1]{LupercioUribe2} that the loop
groupoid  is Morita equivalent to the groupoid
$$\Loop \Xx=
\begin{array}{c}
\left( \bigsqcup_g \PP_g \right) \times G\\
\twodownarrows \\
\left( \bigsqcup_g \PP_g \right)
\end{array}$$
where $\PP_g:= \{\phi: [0,1] \to M \ \ | \ \ \phi(0) g = \phi(1)
\}$ and $G$ acts on the paths in the natural way, i.e. $\{\phi
\cdot k\} (t) = \phi(t) k$  with $\{\phi\cdot k\}(0) k^{-1}gk
=\{\phi \cdot k\}(1)$.

For an arrow $\Lambda = (\phi, g, k)$ in $\Loop \Xx$ between the
paths $(\phi, g)$ and $(\phi \cdot k, k^{-1}gk)$ with $x =
\phi(0)$, the morphism of groupoids $F$  becomes:
\begin{eqnarray*}
F: \Loop \Xx & \to & \complex^\times\\
(\phi, g,k) & \mapsto & \frac{h\left((x,g),(x
g,k)\right)}{h\left((x, k)(x k, k^{-1}gk)\right)}=
\frac{\varepsilon(g,k)}{\varepsilon(k,k^{-1}gk)}
\end{eqnarray*}
and the connection $\Delta$ is equal to zero.

 In this way we obtain a flat line bundle over the loop groupoid $\Loop \Xx$ that once restricted to the inertia groupoid
produces the discrete torsion. This localization procedure is explained in
the next section.

\section{Localization at the fixed points} \label{section Localization at the fixed points}

In \cite{LupercioUribe2} we argued that the inertia groupoid can
be understood as the fixed point set of the action of the real
numbers over the loops (the action shifts the paths by a real
number). This groupoid has for objects the constant paths i.e.
maps $\psi : \real \to \Gg_0$ with $\psi(t) = x \in \Gg_0$ for all
$t$, and for morphisms constant arrows $\Lambda : \real \to
\Gg_1$. This description is equivalent to the one given in the
definition \ref{def.inertiagroupoid}.

As we need to remember the source and the target of the morphisms,
the elements of $(\wedge \Gg)_1$ will be pairs $(v, \alpha) \in
\Gg_2$ such that $v \in \wedge \Gg_0$, $s(v,\alpha)=v$ and $
t(v,\alpha) = \alpha^{-1} v \alpha$.
 The structure maps of $\wedge \Gg$ will be written with the letters $s,t,e,i,m$ to differentiate them from
the ones of $\Gg$.

Hence we have an inclusion of groupoids $j :\wedge \Gg \to \Loop
\Gg$  and so we can pull back the line bundle $(F, \Delta)$
previously described to obtain a line bundle with connection over
the inertia groupoid.

\begin{lemma}
The line bundle $(f,\omega) :=j^*(F,\Delta) = (F,\Delta) |_{\wedge
\Gg}$ over the inertia groupoid $\wedge \Gg$ is flat.
\end{lemma}

\begin{proof}
As the paths representing $\wedge \Gg$ are constant, is easy to
see that
$$f = F|_{\wedge \Gg} \ \ \ \ \ \mbox{and} \ \ \ \ \  \omega = \Delta |_{\wedge \Gg_0} = A|_{\wedge \Gg_0}.$$
From equation \ref{2-cocycle B} we see that $d\omega =0$ because
the maps $\source$ and $\target$ in $\wedge \Gg_0$ are equal. Then
the connection $\omega$ over $\wedge \Gg$ is flat.
\end{proof}

These line bundles are the  representatives on what Ruan has
coined ``inner local systems'' (see \cite{Ruan1}) which he uses to
twist the Chen-Ruan cohomology of orbifolds. In fact, all the
constructions he has of ``inner local systems'' could be done
using the procedure outlined in this paper. We believe the only
relevant local systems are the ones obtained via transgression
from a gerbe with connection.

\begin{definition}
In our terminology an ``inner local system'' is a flat line bundle
$\LL$ over the inertia groupoid $\wedge \Gg$ such that:
\begin{itemize}
\item $\LL$ is trivial once restricted to $\ident(\Gg_0) \subset \wedge \Gg_1$ (i.e. $\LL|_{\ident(\Gg_0)} = 1)$ and
\item $i^* \LL = \LL^{-1}$ where $i : \wedge \Gg \to \wedge \Gg$ is the inverse map (i.e. $(i(v,\alpha) = (\alpha^{-1} v \alpha, \alpha^{-1})$).
\end{itemize}
There is an extra condition in Ruan's definition that is trivially
fulfilled by $\LL$. It just says that if $f :\wedge \Gg \to
\complex^\times$ is the map that contains the information on
transition functions, then $f(\alpha_1)f(\alpha_2) =
f(\alpha_1\alpha_2)$ for composable morphisms; this is true
because $f$ is a morphism of groupoids.
\end{definition}

\begin{proposition}
The line bundle $(f,\omega)$ over $\wedge \Gg$ is an inner local
system for $\Gg$.
\end{proposition}
\begin{proof}
As the paths of $\wedge \Gg$ are constant, from the equation
\ref{definitionF} we see that
$$f(v, \alpha) = \frac{h(v,\alpha)}{h(\alpha, \alpha^{-1} v \alpha)} \ \ \ \  \mbox{for} \ \ \ (v, \alpha) \in \wedge \Gg_1.$$

If $v = \ident(x)$ and $\alpha$ goes from $x$ to $y$ then $f(\ident(x),\alpha) = f (\ident(x),\ident(x))$ and
$f(\alpha,\ident(y))= f(\ident(y),\ident(y))$; this follows from the cocycle condition of $f$ applied to the triples $(\ident(x), \ident(x), \alpha)$
and $(\alpha, \ident(y), \ident(y))$.
Hence $f(v, \alpha) = \frac{h(\ident(x),\ident(x))}{h(\ident(y), \ident(y))}$, which means that the value of the gluing functions
do not depend on the arrow but on its end points. This implies that the restriction of $\LL$ to $\ident(\Gg_0)$ is trivial.

Now as $f$ is a morphism of groupoids, then $f(\Lambda)^{-1} =
f(i\Lambda)$ and hence the second condition holds.
\end{proof}
\subsection{Global quotients}

Recall that for the orbifold $\Xx :=[M/G]$ the inertia groupoid
$\wedge \Xx$ is Morita equivalent to $\sqcup_{(g)} [M^g / C(g)]$
where $M^g$ are the fixed point set of $g$, $C(g)$ is the
centralizer of $g$ in $G$ and the disjoint union runs over $(g)$
the conjugacy classes of elements in $G$.

If we forget the connective structure, the construction outlined
in this paper assigns to every gerbe over $\Xx$ a line bundle over
$\wedge \Xx$. Via this transgression we get $C(g)$ equivariant
line bundles $\LL_g$ over $M^g$.

$$H^3_G(M,\integer) \cong H^3(\Xx, \integer) \longrightarrow H^2(\wedge \Xx, \integer) \cong \bigoplus_{(g)} H^2_{C(g)}(M^g, \integer)$$
 These line bundles
form an inner local system in the sense of Ruan, but also they are
the coefficients Freed-Hopkins-Teleman \cite{FreedHopkinsTeleman}
used to twist the cohomology of the twisted sectors in order to
get a Chern character isomorphism with the twisted K-theory of the
orbifold.

In the case of a gerbe coming from discrete torsion we obtain
$\U{1}$ representations of the groups $C(g)$. These representations
were used by Adem and Ruan \cite{AdemRuan} to twist the orbifold
cohomology and in this way they obtained an isomorphism with the twisted
orbifold $K$-theory \cite{LupercioUribe1}.

\section{Generalized holonomy}

In this last section we want to emphasize that the holonomy map
for gerbes over a groupoid can be generalized to $n$-gerbes.

\begin{theorem}
There is a natural homomorphism
$$\tau_{n} :\breve{C}^n(\Gg,\complex^\times(n+1)_\Gg) \To \breve{C}^{n-1}(\Loop \Gg,\complex^\times(n)_{\Loop \Gg})$$
that sends $n$-cocycles to $(n-1)$-cocycles (i.e. $(n-1)$ gerbes with connection over $\Gg$ to $(n-2)$ gerbes with connection over the loop groupoid),
  commutes with the coboundary operator (i.e. $\tau_n \circ (\delta +(-1)^{n} d) =(\delta +(-1)^{n-1} d ) \circ \tau_{n-1}$)
and therefore induces a map in cohomology
$$\hyper^n\left(\Gg, \complex^\times(n+1)_\Gg \right)
\To \hyper^{n-1}\left(\Loop \Gg, \complex^\times(n)_{\Loop \Gg} \right).$$
\end{theorem}

\begin{proof} Let's define first the map $\tau_n$.
Take $(\omega, \theta^1, \dots, \theta^n) \in \breve{C}^n(\Gg,\complex^\times(n+1)_\Gg)$ with $\omega : \Gg_n \to \complex^\times$
and $\theta^j \in \Gamma(\Gg_{n-j},\Aa^j_{\Gg,\complex})$ and let
$$(F, \Delta^1, \dots , \Delta^{n-1}) :=\tau_n (\omega, \theta^1, \dots, \theta^n)$$
with $F : \Loop \Gg_{n-1} \to \complex^\times$ and  $\Delta^j \in \Gamma(\Loop \Gg_{n-1-j},\Aa^j_{\Loop \Gg,\complex})$ defined in the
following way.

For ${\bf \Lambda} = (\Lambda^1, \dots, \Lambda^{n-1})$ a set of
$n-1$ composable morphisms in $\Loop \Gg_{n-1}$ joining the objects
$\psi^0 , \dots, \psi^{n-1}$ with ${\bf \Lambda}_i : I_i=[\alpha_{i-1},\alpha_i] \to
\Gg_{n-1}$ for $0=\alpha_0 < \alpha_1 <\cdots < \alpha_p =1$, we define
\begin{eqnarray*}
F({\bf \Lambda}) & := & \exp \left( \sum_{i=1}^p \int_{I_i} ({\bf \Lambda}_i)^* \theta^1 \right)\times  \\
& & \prod_{i=1}^p \prod_{j=0}^{n-1} \left(
\omega(\Lambda^1_i(\alpha_i),\dots,\Lambda^j_i(\alpha_i), \psi^j(\alpha_i), \Lambda^{j+1}_{i+1}(\alpha_i), \dots, \Lambda^{n-1}_{i+1}(\alpha_i)) \right)^{{(-1)}^{j+n}}
\end{eqnarray*}

Now let ${\bf \Xi^a} = ( \Xi^{a,1}, \dots, \Xi^{a, n-1-k})$, ${\bf a} \in \{{\bf 1}, \dots,{\bf   k}\}$ 
be vector fields over ${\bf \Lambda} = (\Lambda^1, \dots, \Lambda^{n-1-k}) \in \Loop \Gg_{n-1-k}$ 
 with 
${\bf \Xi^a}_i : I_i \to (T\Gg)_{n-1-k}$,
joining the objects $\xi^{a,0}, \dots, \xi^{a,n-k-1}$ of the tangent loop groupoid; i.e. $\xi^{a,j}$ is a vector
field over $\psi^j$ and $\Xi^{a,j}$ is an arrow between $\xi^{a,j-1}$ and $ \xi^{a,j}$
as well as a vector field over $\Lambda^j$.

For $m \in \{0, \dots, n-1-k\}$ we construct the following set of arrows in $(T \Gg)_{n-k}$:
$$\vartheta_m {\bf \Xi^a} (\alpha_i) := \left( {\Xi}^{a,1}_i(\alpha_i), \dots,  {\Xi}^{a,m}_i(\alpha_i),\xi^m(\alpha_i), 
{\Xi}^{a,m+1}_{i+1}(\alpha_i), \dots , {\Xi}^{a,n-1-k}_{i+1}(\alpha_i) \right).$$
 
Define,
\begin{eqnarray*}
\langle \Delta^k_{\bf \Lambda}, ({\bf \Xi^1}, \dots , {\bf \Xi^k}) \rangle & := &
\sum_{i=1}^p \int_{I_i} \theta^{k+1} \left( \frac{d {\bf \Lambda}_i}{dt}, {\bf \Xi^1}_i(t), \dots , {\bf \Xi^k}_i(t) \right)dt \\
& &+ \sum_{i=1}^p \sum_{m=0}^{n-1-k} (-1)^{m+n}
\left\langle \theta^k  , \left( \vartheta_m {\bf \Xi^1}(\alpha_i), \dots , \vartheta_m{\bf \Xi^k}(\alpha_i)\right)
\right\rangle
\end{eqnarray*}

In what follows we will only show that the map $\tau_n$ sends cocycles to cocycles. The other part of the proof can be done
following the steps of the theorem \ref{theoremlinebundle}.
Let's suppose that $\delta \theta^{k+1} = (-1)^n d \theta^k$ and we want to prove that $\delta \Delta^{k+1} = (-1)^{n-1} d \Delta^{k}$.

Both $\delta \Delta^{k+1}$ and  $d \Delta^{k}$ are in $\breve{C}(\Loop \Gg_{n-1-k}, \Aa^{k+1}_{\Loop \Gg, \complex})$, so we need to
take ${\bf \Xi^a}$ ${\bf a} \in\{ {\bf 1},\dots {\bf k+1}\}$ vector fields over $\Lambda$.

The proof of
$$\left\langle \delta \Delta^{k+1} , \left( {\bf \Xi^1}, \dots, {\bf \Xi^{k+1}} \right) \right\rangle = 
(-1)^{n-1} \left\langle d \Delta^k ,    \left( {\bf \Xi^1}, \dots, {\bf \Xi^{k+1}} \right) \right\rangle$$
will be done by thickening $\Lambda$ in the directions of the ${\bf \Xi^a}$ and then integrating over this tubular neighborhood.

Then let ${\bf \Lambda}( {\vec s} ) \in \Loop \Gg_{n-1-k}$  with ${\vec s} := (s_1, \dots, s_{k+1})$ be such that
${\bf \Lambda} ({\vec 0})= {\bf \Lambda}$ and $\frac{d {\bf \Lambda }(\vec s)}{ds_a}|_{\vec s =0} = {\bf \Xi^a}$. We argue
that
\begin{eqnarray}
\int_{[-\epsilon, \epsilon]^{k+1 } }
\left\langle \delta \Delta^{k+1} ,\left(\frac{d {\bf \Lambda }(\vec s)}{ds_1} , \dots,\frac{d {\bf \Lambda }(\vec s)}{ds_{k+1}}
 \right) \right\rangle d \vec s  = \ \ \ \ \ \ \ \ \ \ \ \ \ \ \  \ \ \ \ \ \ \  \nonumber \\
\ \ \ \ \ \ \ \ \ \ \ \ \ \ \ \  (-1)^{n-1} \int_{[-\epsilon, \epsilon]^{k+1 } }
 \left\langle d \Delta^k ,    \left( \frac{d {\bf \Lambda }(\vec s)}{ds_1} , \dots,\frac{d {\bf \Lambda }(\vec s)}{ds_{k+1}}
 \right) \right\rangle d \vec s. \label{integrals}
\end{eqnarray}
We just need  one last piece of information, the face maps associated to the coboundary operator $\delta$.
They are $\varrho_l : \Gg_{n-k} \to \Gg_{n-k-1}$ with $\varrho_l(g_1, \dots, g_{n-k}) = (g_1, \dots, g_lg_{l+1}, \dots, g_{n-k})$,
$\varrho_0(g_1, \dots, g_{n-k})= (g_2, \dots, g_{n-k})$, $\varrho_{n-k}(g_1, \dots, g_{n-k})= (g_1, \dots, g_{n-k-1})$
and $\rho_l:\Loop \Gg_{n-k-1} \to \Loop \Gg_{n-k-2}$ defined in the same way.

It is easy to see that
\begin{eqnarray*}
\varrho_l \left( \vartheta_m \frac{d {\bf \Lambda}(\vec s)}{ds_a}(\alpha_i) \right)
 & = & \vartheta_{m-1} \left( \rho_l \frac{d {\bf \Lambda}(\vec s)}{ds_a} \right)(\alpha_i) \ \ \ \ \ \mbox{ for $l < m$}\\
\varrho_{m-1} \left( \vartheta_m \frac{d {\bf \Lambda}(\vec s)}{ds_a}(\alpha_i) \right)
 & = & \varrho_{m-1} \left (\vartheta_{m-1} \frac{d {\bf \Lambda}(\vec s)}{ds_a}(\alpha_i) \right) \\
\varrho_l \left( \vartheta_m \frac{d {\bf \Lambda}(\vec s)}{ds_a}(\alpha_i) \right)
 & = & \vartheta_{m} \left( \rho_{l-1} \frac{d {\bf \Lambda}(\vec s)}{ds_a} \right)(\alpha_i) \ \ \ \ \ \mbox{ for $l > m+1$}\\
\end{eqnarray*}
and note that the only elements not paired are 
$$ \varrho_0 \left( \vartheta_0 \frac{d {\bf \Lambda}(\vec s)}{ds_a}(\alpha_i) \right) \ \ \ \mbox{and} \ \ \ \ 
\varrho_{n-k} \left( \vartheta_{n-k-1} \frac{d {\bf \Lambda}(\vec s)}{ds_a}(\alpha_i) \right),$$
these will play an important role in what follows.

Writing ${\bf \Lambda} := {\bf \Lambda}(\vec s)$ we have that the left hand side of \ref{integrals}  becomes:

$$
LHS (\ref{integrals})  =
\int_{[-\epsilon, \epsilon]^{k+1 } } \left (\sum_{i=1}^p \int_{I_i}  (-1)^n d \theta^{k+1} 
\left( \frac{d {\bf \Lambda}_i}{dt},  \frac{d {\bf \Lambda}_i}{ds_1}, \dots,  \frac{d {\bf \Lambda}_i}{ds_{k+1}}
 \right)  dt \right. $$
$$ +  \left.
\sum_{i=1}^p \sum_{m=0}^{n-k-2} \sum_{l=0}^{n-k-1} (-1)^{m+l+n} \left\langle \theta^{k+1}; \vartheta_m \left( \rho_l 
 \frac{d {\bf \Lambda}_i}{ds_1} \right) (\alpha_i), \dots , \vartheta_m \left( \rho_l
\frac{d {\bf \Lambda}_i}{ds_{k+1}} \right) (\alpha_i) \right\rangle \right) d \vec s
$$
after replacing $\delta \theta^{k+2}$ by $(-1)^{n} d \theta^{k+1}$, and composing by  the maps $\rho_l$
from the definition of $\delta$. 
  
And the right hand side becomes:

$$
RHS (\ref{integrals}) = \sum_{j=1}^{k+1} (-1)^{j+n} \sum_{i=1}^p \int_{[-\epsilon, \epsilon]^k}  \int_{I_i}
\ \ \ \ \ \ \ \ \ \ \ \ \ \ \ \ \ \ \ \ \ \ \ \ \ \ $$
$$\ \ \ \ \ \ \ \ \ \ \ \ \ \ \ \ \ \ \ 
\left. \theta^{k+1} \left(   \frac{d {\bf \Lambda}_i}{dt},  \frac{d {\bf \Lambda}_i}{ds_1},
\dots, {\widehat{  \frac{d {\bf \Lambda}_i}{ds_{j}}}}, \dots
 \dots,  \frac{d {\bf \Lambda}_i}{ds_{k+1}}
 \right) \right\vert^{s_j=\epsilon}_{s_j=-\epsilon} dt d\vec s $$
$$+\int_{[-\epsilon, \epsilon]^{k+1}} \sum_{i=1}^p \sum_{m=0}^{n-1-k} \sum_{l=0}^{n-k} 
(-1)^{m+l+3n-1} \ \ \ \ \ \ \ \ \ \ \ \ \ \ \ \ \ \ \ \ \ \ \ \ \ \ \ \ \ \ \ \ $$
$$\ \ \ \ \ \ \  \ \ \ \ \ \ \ \ \ \ \ \ \ \ \left\langle \theta^{k+1};
\varrho_l \left( \vartheta_m \frac{d {\bf \Lambda}}{ds_1}(\alpha_i) \right), \dots ,
\varrho_l \left( \vartheta_m \frac{d {\bf \Lambda}}{ds_{k+1}}(\alpha_i) \right) \right\rangle d \vec s
$$

after evaluating in the boundary of $[-\epsilon,\epsilon]^{k+1}$ for the first
summand, and after replacing $d \theta ^{k}$ by $(-1)^n \delta \theta^{k+1}$ and evaluating it via the maps $\varrho_l$
in the second summand.

Applying Stokes theorem to the first summand of LHS we see that it matches the first summand of RHS except by the
term 
$$(-1)^n \sum_{i=0}^p \int_{[-\epsilon, \epsilon]^{k+1}} \theta^{k+1} 
\left. \left(   \frac{d {\bf \Lambda}_i}{ds_1}, \dots,  \frac{d {\bf \Lambda}_i}{ds_{k+1}}
 \right) \right\vert^{t=\alpha_i}_{t=\alpha_{i-1}} d \vec s.$$  
The second summand of RHS matches the second summand of LHS except by the terms
$$(-1)^{m+l+3n-1}\sum_{i=0}^p\int_{[-\epsilon, \epsilon]^{k+1}} \left\langle \theta^{k+1};
\varrho_l \left( \vartheta_m \frac{d {\bf \Lambda}}{ds_1}(\alpha_i) \right), \dots ,
\varrho_l \left( \vartheta_m \frac{d {\bf \Lambda}}{ds_{k+1}}(\alpha_i) \right) \right\rangle d\vec s$$
when $l=0$, $m=0$ and $l=n-k$, $m=n-k-1$. 
It is not difficult to see now that these last two formulas match. Hence proving that if the tuple
$(\omega, \theta^1, \dots, \theta^n)$
is a cocycle it implies that $(F, \Delta^1, \dots , \Delta^{n-1})$
is also a cocycle. 
\end{proof}

Then we can conclude with the following statement:

\begin{theorem}
There is a natural cochain map $\tau$ of degree -1 (the transgression map)
 $$\tau :\breve{C}^*(\Gg,\complex^\times(n+1)_\Gg) \To \breve{C}^{*-1}(\Loop \Gg,\complex^\times(n)_{\Loop \Gg})$$
that for $*=n$ sends gerbes to gerbes, and induces a map in Deligne cohomology
$$\hyper^*\left(\Gg, \complex^\times(n+1)_\Gg \right)
\To \hyper^{*-1}\left(\Loop \Gg, \complex^\times(n)_{\Loop \Gg} \right).$$
\end{theorem}

Note that for $* \neq n$ we get the topological transgression map (see Prop. \ref{uninteresting part}), i.e.
\begin{eqnarray*}
H^*(\Gg, \integer) &\To& H^{*-1}(\Loop \Gg, \integer)  \ \ \mbox{ for } \  *>n\\
H^*(\Gg, \real/\integer) & \To & H^{*-1}(\Loop \Gg, \complex^\times)  \ \ \mbox{ for } \ *<n
\end{eqnarray*}

\bibliographystyle{gtart}
\bibliography{holonomy}
\end{document}